\magnification=\magstep1
\hsize=16.5 true cm 
\vsize=23.6 true cm
\font\bff=cmbx10 scaled \magstep1

\font\bffg=cmbx10 scaled \magstep3

\font\smc=cmcsc10 
\parindent0cm
%\overfullrule=0cm
%\nopagenumbers
%%%%%%%%%%%%%%%%%%%%%%%%%%%%%%%%
%%%%ABKšRZUNGEN im Textsatz %%%%
           %
            %
\def\bp{\bigskip}              %
\def\mp{\medskip}              %
\def\sp{\smallskip}            %
               %
%%%%%%%%%%%%%%%%%%%%%%%%%%%%%%%%
%%%ZAHLENMENGEN%%%%%%%%%%%%%%%%%
           %
\def\Bbb#1{\hbox{\boldmas #1}} %
\def\R{\Bbb R}                 %
\def\Z{\Bbb Z}                 %
%%%%%%%%%%%%%%%%%%%%%%%%%%%%%%%%%%%%%%%%%%%%%%%%%%%%%%%%%%%%%%%%%%%%%%%
\expandafter\edef\csname amssym.def\endcsname{%
       \catcode`\noexpand\@=\the\catcode`\@\space}

\catcode`\@=11

\def\undefine#1{\let#1\undefined}
\def\newsymbol#1#2#3#4#5{\let\next@\relax
 \ifnum#2=\@ne\let\next@\msafam@\else
 \ifnum#2=\tw@\let\next@\msbfam@\fi\fi
 \mathchardef#1="#3\next@#4#5}
\def\mathhexbox@#1#2#3{\relax
 \ifmmode\mathpalette{}{\m@th\mathchar"#1#2#3}%
 \else\leavevmode\hbox{$\m@th\mathchar"#1#2#3$}\fi}
\def\hexnumber@#1{\ifcase#1 0\or 1\or 2\or 3\or 4\or 5\or 6\or 7\or 8\or
 9\or A\or B\or C\or D\or E\or F\fi}

\font\tenmsa=msam10
\font\sevenmsa=msam7
\font\fivemsa=msam5
\newfam\msafam
\textfont\msafam=\tenmsa
\scriptfont\msafam=\sevenmsa
\scriptscriptfont\msafam=\fivemsa
\edef\msafam@{\hexnumber@\msafam}
\mathchardef\dabar@"0\msafam@39
\def\dashrightarrow{\mathrel{\dabar@\dabar@\mathchar"0\msafam@4B}}
\def\dashleftarrow{\mathrel{\mathchar"0\msafam@4C\dabar@\dabar@}}

\def\ulcorner{\delimiter"4\msafam@70\msafam@70 }
\def\urcorner{\delimiter"5\msafam@71\msafam@71 }
\def\llcorner{\delimiter"4\msafam@78\msafam@78 }
\def\lrcorner{\delimiter"5\msafam@79\msafam@79 }
\def\yen{{\mathhexbox@\msafam@55}}
\def\checkmark{{\mathhexbox@\msafam@58}}
\def\circledR{{\mathhexbox@\msafam@72}}
\def\maltese{{\mathhexbox@\msafam@7A}}

\font\tenmsb=msbm10
\font\sevenmsb=msbm7
\font\fivemsb=msbm5
\newfam\msbfam
\textfont\msbfam=\tenmsb
\scriptfont\msbfam=\sevenmsb
\scriptscriptfont\msbfam=\fivemsb
\edef\msbfam@{\hexnumber@\msbfam}
\def\Bbb#1{{\fam\msbfam\relax#1}}
\def\widehat#1{\setbox\z@\hbox{$\m@th#1$}%
 \ifdim\wd\z@>\tw@ em\mathaccent"0\msbfam@5B{#1}%
 \else\mathaccent"0362{#1}\fi}

\def\widetilde#1{\setbox\z@\hbox{$\m@th#1$}%
 \ifdim\wd\z@>\tw@ em\mathaccent"0\msbfam@5D{#1}%
 \else\mathaccent"0365{#1}\fi}
\font\teneufm=eufm10
\font\seveneufm=eufm7
\font\fiveeufm=eufm5
\newfam\eufmfam
\textfont\eufmfam=\teneufm
\scriptfont\eufmfam=\seveneufm
\scriptscriptfont\eufmfam=\fiveeufm

\newsymbol\risingdotseq 133A
\newsymbol\fallingdotseq 133B
\newsymbol\complement 107B
\newsymbol\nmid 232D
\newsymbol\rtimes 226F
\newsymbol\thicksim 2373

\font\eightmsb=msbm8   \font\sixmsb=msbm6   \font\fivemsb=msbm5
\font\eighteufm=eufm8  \font\sixeufm=eufm6  \font\fiveeufm=eufm5
\font\eightrm=cmr8     \font\sixrm=cmr6     \font\fiverm=cmr5
\font\eightbf=cmbx8    \font\sixbf=cmbx6    
      \font\eighti=cmmi8   \font\sixi=cmmi6
\font\ninesy=cmsy9     \font\eightsy=cmsy8  \font\sixsy=cmsy6
     \font\eightit=cmti8  
     \font\eightsl=cmsl8  
     \font\eighttt=cmtt8
 %SLANTED TYPEWRITER 10 POINT

\font\eightsmc=cmcsc8
\newskip\ttglue
\newfam\smcfam
\def\eightpoint{\def\rm{\fam0\eightrm}%
  \textfont0=\eightrm \scriptfont0=\sixrm \scriptscriptfont0=\fiverm
  \textfont1=\eighti \scriptfont1=\sixi \scriptscriptfont1=\fivei
  \textfont2=\eightsy \scriptfont2=\sixsy \scriptscriptfont2=\fivesy
  \textfont3=\tenex \scriptfont3=\tenex \scriptscriptfont3=\tenex
  \def\smc{\fam\smcfam\eightsmc}
  \textfont\smcfam=\eightsmc          
    %\scriptfont\smcfam=\sixsmc   \scriptscriptfont\smcfam=\fivesmc
\textfont\eufmfam=\eighteufm              \scriptfont\eufmfam=\sixeufm
     \scriptscriptfont\eufmfam=\fiveeufm
\textfont\msbfam=\eightmsb            \scriptfont\msbfam=\sixmsb
     \scriptscriptfont\msbfam=\fivemsb
\def\it{\fam\itfam\eightit}%
  \textfont\itfam=\eightit
  \def\sl{\fam\slfam\eightsl}%
  \textfont\slfam=\eightsl
  \def\bf{\fam\bffam\eightbf}%
  \textfont\bffam=\eightbf \scriptfont\bffam=\sixbf
   \scriptscriptfont\bffam=\fivebf
  \def\tt{\fam\ttfam\eighttt}%
  \textfont\ttfam=\eighttt
  \tt \ttglue=.5em plus.25em minus.15em
  \normalbaselineskip=9pt
  \def\MF{{\manual opqr}\-{\manual stuq}}%
  \let\big=\eightbig
  \setbox\strutbox=\hbox{\vrule height7pt depth2pt width\z@}%
  \normalbaselines\rm}
\def\eightbig#1{{\hbox{$\textfont0=\ninerm\textfont2=\ninesy
  \left#1\vbox to6.5pt{}\right.\n@space$}}}

\catcode`@=13 % reinstating the catcode of @
%%%%%%%%%%%%%%%%%%%%%%%%%%%%%%%%%%%%%%%%%%%%%%%%%%%%%%%%%
%%%%%%%%%%%%%%%%%%%%%%%
\centerline{\bffg Many coarse topologies on the real line}
\medskip\smallskip
\centerline{\bff Gerald Kuba}
\bigskip
\vbox{\eightpoint

{\bf Abstract.} 
Let $\,c=|{\Bbb R}|\,$ denote the cardinality of the continuum
and let $\,\eta\,$ denote the Euclidean topology on $\,{\Bbb R}\,$.
Let $\,{\cal L}\,$ denote the family of all Hausdorff topologies $\,\tau\,$
on $\,{\Bbb R}\,$ with $\,\tau\subset\eta\,$.
Let $\,{\cal L}_1\,$ resp.~$\,{\cal L}_2\,$ resp.~$\,{\cal L}_3\,$ 
denote the family 
of all $\,\tau\in{\cal L}\,$ where $\,({\Bbb R},\tau)\,$ is 
{\it completely normal} resp.~{\it second countable}
resp.~{\it not regular}.
Trivially, $\,{\cal L}_1\cap{\cal L}_3=\emptyset\,$ and
$\,|{\cal L}_i|\leq|{\cal L}|\leq 2^c\,$ and $\,|{\cal L}_2|\leq c\,$.
For $\,\tau\in{\cal L}\,$ the space $\,({\Bbb R},\tau)\,$ is metrizable 
if and only if $\,\tau\in{\cal L}_1\cap{\cal L}_2\,$.
We show that, up to homeomorphism, both
$\,{\cal L}_1\,$ and $\,{\cal L}_3\,$ contain precisely $\,2^c\,$
topologies and 
$\,{\cal L}_2\,$ contains precisely $\,c\,$ completely metrizable topologies.
For $\,2^c\,$ non-homeomorphic
topologies $\,\tau\in{\cal L}_1\,$ the space $\,({\Bbb R},\tau)\,$
is {\it Baire}, but there are also  $\,2^c\,$ non-homeomorphic topologies
$\,\tau\in{\cal L}_1\,$  and $\,c\,$ non-homeomorphic topologies
$\,\tau\in{\cal L}_1\cap{\cal L}_2\,$                
where $\,({\Bbb R},\tau)\,$ is of {\it first category}.
Furthermore, we investigate the {\it complete lattice} $\,{\cal L}_0\,$
of all topologies $\,\tau\in{\cal L}\,$ such that $\,\tau\,$ and $\,\eta\,$ 
coincide on $\,{\Bbb R}\setminus\{0\}\,$. In the lattice $\,{\cal L}_0\,$ 
we find $\,2^c\,$ (non-homeomorphic)
immediate predecessors of the maximum $\,\eta\,$,
whereas the minimum of $\,{\cal L}_0\,$ is a compact topology 
without immediate successors in $\,{\cal L}_0\,$.
We construct chains of homeomorphic topologies
in $\,{\cal L}_0\cap{\cal L}_1\cap{\cal L}_2\,$ and in 
$\,{\cal L}_0\cap{\cal L}_2\cap{\cal L}_3\,$ and in 
$\,{\cal L}_0\cap({\cal L}_1\setminus{\cal L}_2)\,$ and in 
$\,{\cal L}_0\cap({\cal L}_3\setminus{\cal L}_2)\,$ such that
the length of each chain is $\,c\,$ (and hence maximal).
We also track down a chain in $\,{\cal L}_0\,$ of length
$\,2^\lambda\,$ where $\,\lambda\,$ 
is the smallest cardinal number $\,\kappa\,$ with $\,2^\kappa>c\,$.}
\bigskip
{\bff 1. Introduction} 
\medskip
Write $\,|S|\,$ for the cardinality (the size) of a the set $\,S\,$ and let 
$\,c=|{\Bbb R}|\,$ denote the cardinality of the continuum.
Let $\,\eta\,$ denote the Euclidean topology on $\,{\Bbb R}\,$
and let $\,{\cal L}\,$ denote the family of all topologies $\,\tau\,$
on $\,{\Bbb R}\,$ where $\,\tau\,$ is coarser than $\,\eta\,$
(i.e.~$\,\tau\,$ is a subset of $\,\eta\,$) and
$\,({\Bbb R},\tau)\,$ is a Hausdorff space. 
If $\,\tau\in{\cal L}\,$
and $\,B\,$ is a nonempty bounded 
subset of $\,{\Bbb R}\,$ then the relative topologies of $\,\tau\,$
and $\,\eta\,$ coincide on $\,B\,$.
(Because they coincide on the 
interval $\;[\inf B,\sup B]\;$
due to the well-known fact that a topology cannot be T$_2$ if it is 
strictly coarser than a T$_2$-compact topology.)
Nevertheless, on the whole space $\,{\Bbb R}\,$ the two 
topologies $\,\tau\,$ and $\,\eta\,$ need not coincide.
In fact, as we will see, $\,|{\cal L}|=2^c\,$.
(Note that $\,|{\cal L}|\leq 2^c\,$ is trivial because
$\,|\eta|=c\,$.) Moreover, as we will prove in Section 4, 
$\,{\cal L}\,$ contains $\,2^c\,$ mutually non-homeomorphic topologies 
$\,\tau\,$ such that 
$\,({\Bbb R},\tau)\,$ is a completely normal Baire 
space. In Section 8 we will prove that 
$\,{\cal L}\,$ also contains $\,2^c\,$ mutually non-homeomorphic topologies 
$\,\tau\,$ such that 
$\,({\Bbb R},\tau)\,$ is a completely normal space of first category. 
\medskip
For every $\,\tau\in{\cal L}\,$
the space $\,({\Bbb R},\tau)\,$ is separable and arcwise connected
and $\sigma$-compact.
Separability is trivial since $\,{\Bbb Q}\,$ is clearly a dense set
in $\,({\Bbb R},\tau)\,$. Arcwise connectedness and $\sigma$-compactness
follow immediately
from the coincidence of $\,\eta\,$ and $\,\tau\,$ on each Euclidean 
compact interval.
Whereas the Euclidean space $\,{\Bbb R}\,$ is second countable,
for arbitrary $\,\tau\in{\cal L}\,$ the space
$\,({\Bbb R},\tau)\,$ need not be second countable.
In fact, there cannot be more than $\,c\,$ second countable 
topologies in the family $\,{\cal L}\,$ 
since $\,|\eta|=c\,$ and a set of size $\,c\,$ has 
precisely $\,c\,$ countable subsets. 
Due to separability, for $\,\tau\in{\cal L}\,$
the space $\,({\Bbb R},\tau)\,$ is metrizable if and only if 
it is regular and second countable. 
In particular, there are at most $\,c\,$ 
metrizable topologies in the family $\,{\cal L}\,$.
In Section 7 we will prove that 
there exist $\,c\,$ mutually non-homeomorphic 
topologies $\,\tau\in{\cal L}\,$ such that
$\,({\Bbb R},\tau)\,$ is completely metrizable.
In Section 9 we will prove that 
there exist $\,c\,$ mutually non-homeomorphic 
topologies $\,\tau\in{\cal L}\,$ such that
$\,({\Bbb R},\tau)\,$ is a metrizable space of first category.
\medskip
Let us call the
image $\,g({\Bbb R})\,$ of 
any continuous one-to-one mapping $\,g\,$ from the Euclidean space 
$\,{\Bbb R}\,$ 
into a Hausdorff space $\,X\,$ a {\it real arc}.
There is a natural correspondence 
between topologies in the family $\,{\cal L}\,$ 
and real arcs.
Because, with $\,g\,$ and $\,X\,$ as above, evidently
the family $\,\tau_g\,$ of all sets $\,g^{-1}(U)\,$ 
where $\,U\,$ is an open subset of $\,X\,$ is 
a topology in the family $\,{\cal L}\,$
and $\,g\,$ defines a homeomorphism 
between the space $\,({\Bbb R},\tau_g)\,$ and the subspace 
$\,g({\Bbb R})\,$ of $\,X\,$. Conversely,
for each $\,\tau\in{\cal L}\,$ the space 
$\,({\Bbb R},\tau)\,$ is a real arc 
since the identity is a continuous mapping from 
$\,({\Bbb R},\eta)\,$ onto $\,({\Bbb R},\tau)\,$.
As a consequence of our enumeration results mentioned above 
and proved in Sections 4 and 7,
up to homeomorphism there are precisely $\,2^c\,$ 
completely normal real arcs
and precisely $\,c\,$ 
completely metrizable  real arcs.
Our result on completely metrizable topologies in $\,{\cal L}\,$
will be proved by constructing real arcs within
the Euclidean space $\,{\Bbb R}^3\,$.
\bigskip
{\bff 2. Locally and globally coarse topologies}
\medskip
If $\,\tau\,$ is a topology on the set $\,{\Bbb R}\,$ and $\,a\in{\Bbb R}\,$ 
then let $\,{\cal N}_\tau(a)\,$ denote the filter of the neighborhoods
of the point $\,a\,$ in the space $\,({\Bbb R},\tau)\,$. 
Trivially, $\,{\cal N}_\tau(a)\subset{\cal N}_\eta(a)\,$
for every $\,\tau\in{\cal L}\,$.
Let us call a topology $\,\tau\,$ in our family $\,{\cal L}\,$ 
{\it coarse at
the point} $\,a\in{\Bbb R}\,$ if and only if 
$\,{\cal N}_\tau(a)\not={\cal N}_\eta(a)\,$.
A proof of the following lemma is straightforward.
\smallskip
{\bf Lemma 1.} {\it If an injective mapping $\,g\,$ with domain $\,{\Bbb R}\,$
defines a real arc $\,g({\Bbb R})\,$ then  
the topology $\,\tau_g\,$ in $\,{\cal L}\,$  corresponding with $\,g\,$ 
is coarse at $\,a\in {\Bbb R}\,$ if and only if the bijection
$\,g^{-1}\,$ from $\,g({\Bbb R})\,$ 
onto $\,{\Bbb R}\,$ is not continuous at $\,g(a)\,$.}
\smallskip
The following proposition makes 
it easy to detect whether a topology $\,\tau\in{\cal L}\,$ 
is coarse at a point $\,a\in{\Bbb R}\,$.
\smallskip
{\bf Proposition 1.} {\it A topology $\,\tau\in{\cal L}\,$ 
is coarse at a point $\,a\in{\Bbb R}\,$ if and only if
every set in the filter $\,{\cal N}_\tau(a)\,$ 
is an unbounded subset of $\,{\Bbb R}\,$.}
\smallskip
{\it Proof.} Let $\,\tau\in{\cal L}\,$ and $\,a\in{\Bbb R}\,$ 
and assume that some $\,U\in{\cal N}_\tau(a)\,$ is bounded.
Fix $\,\delta>0\,$ so that 
$\,[a-\delta,a+\delta]\subset U\,$
and let $\,0<\varepsilon\leq \delta\,$ be arbitrary.
The Euclidean compact set
$\;A\,=\,[\inf U,a-\varepsilon]\cup[a+\varepsilon, \sup U]\;$
is compact and hence closed in the space $\,({\Bbb R},\tau)\,$.
Consequently, $\;]a-\varepsilon,a+\varepsilon[\;=\,U\setminus A\,$ 
is $\tau$-open whenever $\,0<\varepsilon\leq \delta\,$
and hence $\,{\cal N}_\tau(a)={\cal N}_\eta(a)\,$, {\it q.e.d.}
\medskip
The following proposition provides a nice and very useful 
characterization of the first-category 
topologies in the family $\,{\cal L}\,$.
\smallskip
{\bf Proposition 2.} {\it For $\,\tau\in{\cal L}\,$ the 
space $\,({\Bbb R},\tau)\,$ is of first category if and only if 
every nonempty open set in the space $\,({\Bbb R},\tau)\,$
is an unbounded subset of $\,{\Bbb R}\,$.}
\smallskip
{\it Proof.} Assume 
firstly that $\,\tau\in{\cal L}\,$ and every nonempty $\tau$-open set is 
unbounded.
Then for each $\,n\in{\Bbb N}\,$ 
the set $\,[-n,n]\,$ is  nowhere dense
in the space $\,({\Bbb R},\tau)\,$. (Note that  
the Euclidean compact set $\,[-n,n]\,$ is $\tau$-compact and hence 
$\tau$-closed.) 
Thus the space $\,({\Bbb R},\tau)\,$ is of first category
since $\;{\Bbb R}\,=\,\bigcup_{n=1}^\infty[-n,n]\,$.
%\smallskip
Assume secondly that $\,\tau\in{\cal L}\,$ 
and that $\,({\Bbb R},\tau)\,$ is a space of first category
and suppose that
there would exist a nonempty $\tau$-open set $\,U\,$
which is bounded. As an open subspace of a space of first category,
the set $\,U\,$ equipped with the relative topology of $\,\tau\,$
would be a space of first category. But this space is identical with 
$\,U\,$ equipped with the relative topology of $\,\eta\,$
(since $\,U\,$ is bounded)
and, naturally, the Euclidean space $\,U\,$ is of second category.
This contradiction finishes the proof, {\it q.e.d.}
\medskip
{\it Remark.} As a trivial consequence of Propositions 1 and 2, for 
$\,\tau\in{\cal L}\,$ the space $\,({\Bbb R},\tau)\,$ is of first category
if and only if $\,\tau\,$ is {\it everywhere} coarse.
In [5] we construct $\,2^{2^c}\,$ non-homeomorphic
connected topologies $\,\tau\,$ on $\,{\Bbb R}\,$ with certain properties
where $\,\tau\,$ is finer than $\,\eta\,$.
In [5] it is not explicitly stated that all
these topologies $\,\tau\,$ are actually 
{\it everywhere} finer than $\,\eta\,$, 
i.e.~$\,{\cal N}_\eta(a)\,$ is a proper subset of $\,{\cal N}_\tau(a)\,$
for every $\,a\in {\Bbb R}\,$. However, some of these $\,2^{2^c}\,$ topologies
are of first category, but some of them are of second category.
\medskip
For $\,\tau\in{\cal L}\,$ let $\,C(\tau)\,$ denote the set of all points $\,a\,$ 
such that $\,\tau\,$ is coarse at $\,a\,$. 
Clearly, if $\,C(\tau)\not={\Bbb R}\,$ then the subspace topologies 
of $\,\tau\,$ and $\,\eta\,$ coincide on the set $\,{\Bbb R}\setminus C(\tau)\,$.
The following proposition shows that 
the set $\,C(\tau)\,$ is always of a very special form.
\medskip
{\bf Proposition 3.} {\it Let $\,\tau\in{\cal L}\,$. Then $\,C(\tau)\,$
is a closed subset of the Euclidean space $\,{\Bbb R}\,$.
Moreover, the set $\,C(\tau)\,$ is 
closed and meager in the space $\,({\Bbb R},\tau)\,$.}
\smallskip
{\it Proof.} Let $\,\tau\in{\cal L}\,$. 
Firstly we verify that $\,C(\tau)\,$ is closed
in the space $\,({\Bbb R},\tau)\,$.
(Then, of course, $\,C(\tau)\,$ is closed in the 
Euclidean space automatically.) 
Assume that $\,x\in{\Bbb R}\,$ is 
a limit point of the set $\,C(\tau)\,$ in the space $\,({\Bbb R},\tau)\,$.
Then $\;U\cap C(\tau)\not=\emptyset\;$ 
for every $\tau$-open set $\,U\,$ in the filter $\,{\cal N}_\tau(x)\,$
and hence every set in the filter $\,{\cal N}_\tau(x)\,$ 
lies in the filter $\,{\cal N}_\tau(a)\,$ for some $\,a\in C(\tau)\,$. 
Thus every set in $\,{\cal N}_\tau(x)\,$ is unbounded 
by Proposition~1. Hence $\,x\in C(\tau)\,$ by Proposition 1.
Therefore the set $\,C(\tau)\,$ is $\tau$-closed.
Since $\,[-n,n]\,$ is compact and hence closed 
in the space $\,({\Bbb R},\tau)\,$ for every $\,n\in{\Bbb N}\,$, all sets 
$\,C(\tau)\cap[-n,n]\,$ are closed in the space $\,({\Bbb R},\tau)\,$.
No point in $\,C(\tau)\cap[-n,n]\,$
is an $\tau$-interior point of $\,C(\tau)\cap[-n,n]\,$
because if $\,a\in C(\tau)\,$ then
$\;S\not\subset[-n,n]\;$ for every $\,S\in{\cal N}_\tau(a)\,$
by Proposition~1. Consequently,
$\,C(\tau)\cap[-n,n]\,$ is nowhere dense in the space $\,({\Bbb R},\tau)\,$
for every $\,n\in{\Bbb N}\,$
and hence the set 
$\;C(\tau)\,=\,\bigcup_{n=1}^\infty(C(\tau)\cap[-n,n])\;$
is meager in the space $\,({\Bbb R},\tau)\,$, {\it q.e.d.}
\medskip
The following proposition generalizes the special fact that 
$\,({\Bbb R},\eta)\,$ is a Baire space with $\,C(\eta)=\emptyset\,$
and will be useful for the proof of the enumeration results in 
Sections 4 and 5.
\medskip
{\bf Proposition 4.} {\it If $\,\tau\in{\cal L}\,$ such that $\,C(\tau)\,$
is a meager set in the space $\,({\Bbb R},\eta)\,$ then 
$\,({\Bbb R},\tau)\,$ is a Baire space.}
\smallskip
{\it Proof.} For $\,\tau\in{\cal L}\,$ assume that 
$\,C(\tau)\,$ is a meager subset of Euclidean space $\,{\Bbb R}\,$.
Then $\,C(\tau)\not={\Bbb R}\,$ and hence 
$\;U\,:=\,{\Bbb R}\setminus C(\tau)\;$ is nonempty.
By Proposition~3 the set $\,U\,$ is Euclidean open (even $\tau$-open).
As an open subspace of the Baire space 
$\,({\Bbb R},\eta)\,$, the space $\,(U,\eta)\,$ is Baire.
The spaces $\,(U,\eta)\,$ and $\,(U,\tau)\,$ are identical 
in view of $\,U\cap C(\tau)=\emptyset\,$ and the definition 
of the set $\,C(\tau)\,$. 
In particular, the space $\,(U,\tau)\,$ is Baire.
As the complement of a meager set, $\,U\,$ is dense in 
the Euclidean space $\,{\Bbb R}\,$ and hence
dense in the space $\,({\Bbb R},\tau)\,$ a fortiori.
This is enough in view of the well-known fact (cf.~[2] 3.9.J.b)
that a Hausdorff 
space must be Baire if some dense subspace is Baire, {\it q.e.d.}
\medskip
The following proposition, which implies that $\,{\cal L}\,$ contains
$\,c\,$ completely metrizable topologies, demonstrates
that the converse of Proposition 4 would be far from being true.
\medskip
{\bf Proposition 5.} {\it For every $\,z\in{\Bbb R}\,$ 
there exists a 
topology $\,\tau_z\in{\cal L}\,$ with $\;C(\tau_z)=\,]{-\infty,z}]\;$
such that all spaces $\,({\Bbb R},\tau_z)\,$
are completely metrizable and homeomorphic.}
\smallskip
{\it Proof.} We work with real arcs
and define for every $\,z\in{\Bbb R}\,$ an injective and continuous mapping 
$\,g_z\,$ from the Euclidean space $\,{\Bbb R}\,$ 
into the Euclidean plane $\,{\Bbb R}^2\,$ by putting
$\;g_z(t)=(t,0)\;$ for $\;t\leq z\;$
and $\;g_z(t)=(z+(t-z)(z+1-t),t-z)\;$ for $\;z\leq t\leq z+1\;$
and $\;g_z(t)\,=\,(z+(z+1-t)|\sin(z+1-t)|\,,\,e^{z+1-t})\;$ for 
$\;t\geq z+1\,$. Clearly, $\,g_z({\Bbb R})\,$
is a closed subset of the complete metric space $\,{\Bbb R}^2\,$. 
We observe that $\,g_z^{-1}\,$ is continuous at 
$\,g_z(a)\,$ if and only if $\;a\in\,]z,\infty[\,$.
(Hence $\,C(\tau_z)=\,]{-\infty,z}]\,$ for $\,\tau_z\in{\cal L}\,$
corresponding with $\,g_z\,$.)
Finally, for every $\,z\in{\Bbb R}\,$ the space 
$\,g_z({\Bbb R})\,$ is homeomorphic to the space $\,g_0({\Bbb R})\,$ 
since the translation $\,(x,y)\mapsto (x-z,y)\,$
of the vector space $\,{\Bbb R}^2\,$ 
maps $\,g_z({\Bbb R})\,$ onto $\,g_0({\Bbb R})\,$, {\it q.e.d.}
\bigskip
{\bff 3. Selecting non-homeomorphic topologies}
\medskip
{\bf Lemma 2.} {\it If 
$\,{\cal H}\subset{\cal L}\,$ and all topologies in $\,{\cal H}\,$
are homeomorphic then $\,|{\cal H}|\leq c\,$.}
\smallskip
{\it Proof.} Firstly, if $\,\tau_1,\tau_2\in{\cal L}\,$
then each continuous function from the space $\,({\Bbb R},\tau_1)\,$ into the 
space $\,({\Bbb R},\tau_2)\,$ 
is completely determined by its values at the points in 
the $\tau_1$-dense set $\,{\Bbb Q}\,$.
Secondly, there are precisely $\,c\,$ functions from $\,{\Bbb Q}\,$ into 
$\,{\Bbb R}\,$, {\it q.e.d.}
\medskip
The following lemma 
makes it very easy to provide mutually non-homeomorphic topologies 
in certain situations.
\medskip
{\bf Lemma 3.} {\it If the size of a family $\,{\cal K}\subset{\cal L}\,$
is greater than $\,c\,$ then $\,{\cal K}\,$ contains a
family $\,{\cal K}'\,$ equipollent to $\,{\cal K}\,$
such that all topologies in $\,{\cal K}\,$ are mutually 
non-homeomorphic.}
\smallskip
{\it Proof.} Define an equivalence relation $\,\sim\,$ on $\,{\cal K}\,$
by putting $\;\tau_1\sim\tau_2\;$ for $\,\tau_i\in{\cal K}\,$
when the spaces $\,({\Bbb R},\tau_1)\,$ and $\,({\Bbb R},\tau_2)\,$ are 
homeomorphic. By Lemma 2 the size of an equivalence class cannot exceed 
$\,c\,$. Consequently, from $\,|{\cal K}|>c\,$ we derive that
the total number of the equivalence classes 
must be $\,|{\cal K}|\,$. So we are done by choosing for $\,{\cal K}'\,$
a set of representatives with respect to the equivalence relation 
$\,\sim\,$, {\it q.e.d.}
\bigskip
{\bff 4. Completely normal Baire topologies}
\medskip\smallskip
The following lemma is very useful in order to 
avoid a lengthy verification of complete normality
by verifying regularity only.
\medskip
{\bf Lemma 4.} {\it Let $\,z\in{\Bbb R}\,$ and $\,\tau\in{\cal L}\,$ with 
$\,C(\tau)=\{z\}\,$. Then the space 
$\,({\Bbb R},\tau)\,$ is second countable if and only if 
some local basis at the point $\,z\,$ is countable.
And the space $\,({\Bbb R},\tau)\,$ is completely normal if and only if it is
regular.}
\medskip
{\it Proof.} Clearly, $\,z\not\in V\in\eta\,$ implies $\,V\in\tau\,$.
This settles the first statement and has also the consequence
that $\;U\cup V\in\tau\;$ whenever
$\,z\in U\in\tau\,$ and $\,V\in\eta\,$.
Assume that $\,({\Bbb R},\tau)\,$ is regular and that
in the  space $\,({\Bbb R},\tau)\,$
we have $\;\overline  A\cap B=A\cap\overline B=\emptyset\;$
for $\,A,B\subset{\Bbb R}\,$. If $\;z\not\in A\cup B\;$
then $\,A\,$ and $\,B\,$ can be separated by $\eta$-open subsets 
of $\,{\Bbb R}\setminus\{z\}\,$ which must be $\tau$-open.
So assume $\,z\in A\cup B\,$ and, say, $\,z\in A\,$.
Then we can find 
disjoint sets $\,U_1,V_1\in\eta\,$ with $\,z\not\in U_1\cup V_1\,$
such that $\;A\setminus\{z\}\,\subset\,U_1\;$ and $\;B\subset V_1\,$.
Furthermore, since the space $\,({\Bbb R},\tau)\,$ is regular, 
we can find disjoint sets $\,U_2,V_2\in\tau\,$ with $\,z\in U_2\,$
and $\,\overline B\subset V_2\,$. Then $\,U_1\cup U_2\,$ and
$\,V_1\cap V_2\,$ are disjoint $\tau$-open sets 
and  $\;A\subset U_1\cup U_2\;$ and  $\;B\subset V_1\cap V_2\,$, {\it q.e.d.}
\medskip\smallskip
Our first main result is the following theorem.
\medskip
{\bf Theorem 1.}\quad {\it There exists a family $\;{\cal T}\subset{\cal L}\;$ 
with $\,|{\cal T}|=2^c\,$ such that
$\,({\Bbb R},\tau)\,$ is a completely normal Baire 
space for each $\,\tau\in{\cal T}\,$
and two spaces $\,({\Bbb R},\tau)\,$ and $\,({\Bbb R},\tau')\,$ are
never homeomorphic for distinct topologies 
$\,\tau,\tau'\in{\cal T}\,$.}
\medskip
{\it Proof.} The cardinal number $\,2^c\,$ indicates that 
the natural way to define $\,{\cal T}\,$ is to use ultrafilters
on a countably infinite set.
It is well-known (see [1])
that an infinite set of size $\,\kappa\,$
carries precisely $\,2^{2^{\kappa}}\,$ free ultrafilters.
In particular, there are $\,2^c\,$ free ultrafilters on $\,{\Bbb Z}\,$. 
Note that no free ultrafilter contains a finite set.
\smallskip
For each free ultrafilter $\,{\cal F}\,$ on $\,{\Bbb Z}\,$
define a topology $\,\tau=\tau[{\cal F}]\,$ on $\,{\Bbb R}\,$
by declaring $\,U\subset{\Bbb R}\,$ open if and only if 
$\,U\,$ is Euclidean open and satisfies $\;0\not\in U\;$
or $\;U\cap{\Bbb Z}\,\in\,{\cal F}\,$.
It is plain that $\,\tau\,$
is a well-defined topology on $\,{\Bbb R}\,$ coarser than $\,\eta\,$.
Further, $\,({\Bbb R},\tau)\,$ is a Hausdorff space, 
whence $\,\tau\in{\cal L}\,$, because if $\,u<v\,$ then 
the intersection $\,{\Bbb Z}\setminus[u,v]\,$ of $\,{\Bbb Z}\,$ and the Euclidean open set 
$\,{\Bbb R}\setminus[u,v]\,$ must lie in $\,{\cal F}\,$
(since $\,{\Bbb Z}\cap[u,v]\,$ is a finite set
and the ultrafilter $\,{\cal F}\,$ is free).
By Proposition 1
we have $\,0\in C(\tau)\,$ 
since $\;M\cap{\Bbb Z}\,\in\,{\cal F}\;$ for every $\,M\in{\cal N}_\tau(0)\,$
and every  $\,S\in{\cal F}\,$ is
an infinite set.
Moreover, $\,C(\tau)=\{0\}\,$ since $\,\tau\,$ and $\,\eta\,$ coincide on 
the Euclidean open set $\,{\Bbb R}\setminus\{0\}\,$.
Hence $\,({\Bbb R},\tau)\,$ is a Baire space by Proposition 4.
\smallskip
We claim that $\,({\Bbb R},\tau)\,$ is completely normal.
By Lemma 4 it is enough to check the 
T$_3$-separation property. 
Let $\,A\subset{\Bbb R}\,$ be $\tau$-closed (and hence $\eta$-closed)
and let $\;b\,\in\,{\Bbb R}\setminus A\,$.
If $\,b\not=0\,$ then we can find $\,\epsilon>0\,$
and $\,U\in\eta\,$ disjoint from $\;V:=\,]b-\epsilon,b+\epsilon[\;$
with $\,0\not\in V\,$
and $\,A\subset U\,$. Then $\,V\,$ is $\tau$-open
and $\;U_\epsilon\,:=\,U\cup({\Bbb R}\setminus[b-\epsilon,b+\epsilon])\;$
is $\tau$-open 
and $\;b\in V\;$ and $\;A\subset U_\epsilon\;$ 
and  $\;U_\epsilon\cap V=\emptyset\,$.
(The set $\,U_\epsilon\cap {\Bbb Z}\,$ lies in the free ultrafilter $\,{\cal F}\,$ 
since $\,{\Bbb Z}\setminus U_\epsilon\,$ is finite.)
If $\,b=0\,$ then $\;B\,:=\,\{0\}\cup({\Bbb Z}\setminus A)\;$ 
is $\eta$-closed and disjoint from $\,A\,$ and hence
we can choose disjoint $\eta$-open sets $\,U,V\,$ 
with $\,A\subset U\,$ and $\,b\in B\subset V\,$.
The set $\,U\,$ is $\tau$-open because $\,0\not\in U\,$
since $\,0\in V\,$ and $\,U\cap V=\emptyset\,$.
The set $\,V\,$ is $\tau$-open because 
$\,{\Bbb Z}\setminus A\in{\cal F}\,$ (since $\,A\,$ is $\tau$-closed) 
and hence from $\,V\cap{\Bbb Z}\supset B\cap{\Bbb Z}\supset{\Bbb Z}\setminus A\,$
 we derive $\,V\cap{\Bbb Z}\in{\cal F}\,$.
\smallskip
Finally we observe that              
$\,\tau[{\cal F}_1]\not\subset\tau[{\cal F}_2]\,$
(and hence $\,\tau[{\cal F}_1]\not=\tau[{\cal F}_2]\,$)
whenever $\,{\cal F}_1\,$ and $\,{\cal F}_2\,$
are distinct free ultrafilters on $\,{\Bbb Z}\,$.
Indeed, if $\,{\cal F}_1\,$ and $\,{\cal F}_2\,$ are free ultrafilters 
on $\,{\Bbb Z}\,$ and 
$\,\tau[{\cal F}_1]\subset\tau[{\cal F}_2]\,$
and $\,S\in{\cal F}_1\,$ then the 
$\tau[{\cal F}_1]$-open set 
$\;\,W\,:=\,\big]{-{1\over 3},{1\over 3}}\big[\,\cup\bigcup_{s\in S}
\big]s-{1\over 3},s+{1\over 3}\big[\,\;$
is a $\tau[{\cal F}_2]$-open neighborhood of $\,0\,$ and hence 
$\;S\cup\{0\}\,=\,W\cap{\Bbb Z}\;$ 
lies in $\,{\cal F}_2\,$, whence $\,S\in{\cal F}_2\,$.
(Note that $\,{\Bbb Z}\setminus\{0\}\in{\cal F}_2\,$ since 
the ultrafilter $\,{\cal F}_2\,$ is free.)
Thus $\,{\cal F}_1\subset{\cal F}_2\,$
and hence $\,{\cal F}_1={\cal F}_2\,$ 
since $\,{\cal F}_1\,$ and $\,{\cal F}_2\,$
are ultrafilters, {\it q.e.d.}
\medskip                      
{\it Remark.} Since $\,{\cal L}\,$ contains only $\,c\,$ second countable 
topologies, there are 
$\,2^c\,$ free ultrafilters $\,{\cal F}\,$ on $\,{\Bbb Z}\,$ 
such that the space $\,({\Bbb R},\tau[{\cal F}])\,$ is not second countable
or, equivalently, that any local basis at $\,0\,$ is uncountable.
In fact, this is true for every free ultrafilter $\,{\cal F}\,$ on $\,{\Bbb Z}\,$. 
Indeed, assume on the contrary 
that the countable family $\;\{\,B_1,B_2,B_3,...\,\}\;$
is a local basis                                         
at $\,0\,$ in the space $\,({\Bbb R},\tau[{\cal F}])\,$.
Then we may choose a sequence $\;a_1,a_2,a_3,...\,$ of distinct
integers and $\;0<\epsilon_n<{1\over 3}\,(n\in{\Bbb N})\;$
such that $\;a_n\,\in\,B_n\setminus\{a_k\,|\,k<n\}\;$
and $\;[a_n-\epsilon_n,a_n+\epsilon_n]\subset B_n\;$
for every $\,n\in{\Bbb N}\,$. 
Then with $\;S\,=\,{\Bbb Z}\setminus\{a_1,a_2,a_3,...\,\}\;$
the set 

\centerline{$\;U\;:=\;\bigcup\limits_{n=1}^\infty]a_n-\epsilon_n,a_n+\epsilon_n[
\;\,\cup\,\bigcup\limits_{s\in S}]s-{1\over 3},s+{1\over 3}[\;$}

is a $\tau[{\cal F}]$-open $\tau[{\cal F}]$-neighborhood 
of $\,0\,$ (since $\;U\cap{\Bbb Z}={\Bbb Z}\in{\cal F}\,$)
with $\;a_n+\epsilon_n\,\in\,B_n\setminus U\;$ and hence
$\;B_n\not\subset U\;$ for every $\,n\in{\Bbb N}\,$.
Thus $\;\{\,B_1,B_2,B_3,...\,\}\;$
is not a local basis at $\,0\,$.
\bigskip
{\bff 5. Non-regular Baire topologies}
\medskip
In view of Theorem 1 and Lemma 4 there arises the question 
whether $\,{\cal L}\,$ contains also $\,2^c\,$ topologies $\,\tau\,$ 
which are Baire because of $\,C(\tau)=\{0\}\,$ and where $\,({\Bbb R},\tau)\,$
is not regular. This is indeed true.
\medskip
{\bf Theorem 2.} {\it There exist $\,2^c\,$ 
mutually non-homeomorphic 
topologies $\,\tau\in{\cal L}\,$  
such that $\,({\Bbb R},\tau)\,$ is a Baire space
which is not regular.}
\medskip
{\it Proof.} It is enough to modify the proof of Theorem 1 
in the following way.
For any free ultrafilter $\,{\cal F}\,$ on $\,{\Bbb Z}\,$
define a topology $\,\sigma[{\cal F}]\,$ on $\,{\Bbb R}\,$
by declaring $\,U\subset{\Bbb R}\,$ open if and only if 
$\,U\,$ is Euclidean open and $\;0\not\in U\;$
or $U\;\supset\;\bigcup_{s\in S}
\big]s-{1\over 3},s+{1\over 3}\big[\;$
for some $\,S\in{\cal F}\,$.
Certainly, 
$\,\sigma[{\cal F}]\,$ is well-defined and Hausdorff.
The space $\,({\Bbb R},\sigma[{\cal F}])\,$ is not regular since, 
for example, the point $\,0\,$ and the obviously
$\sigma[{\cal F}]$-closed set 
$\;\bigcup_{k=-\infty}^{k=\infty}
\big[k+{1\over 3},k+{2\over 3}\big]\;$
cannot be separated by $\sigma[{\cal F}]$-open sets.
Finally, similarly as in the proof of Theorem 1,
$\;\sigma[{\cal F}]\not=\sigma[{\cal F}']\;$
whenever $\,{\cal F}\,$ and $\,{\cal F}'\,$ are distinct
free ultrafilters on $\,{\Bbb Z}\,$, {\it q.e.d.}
\bigskip
{\it Remark.} In the proof of Theorem 1 or Theorem 2 one cannot avoid
an application of Lemma 3 (or a similar transfinite counting argument).
Actually, for every free ultrafilter $\,{\cal F}_0\,$ on $\,{\Bbb Z}\,$
there is an infinite family $\,{\cal U}\,$ of free
ultrafilters  on $\,{\Bbb Z}\,$ with $\,{\cal F}_0\in{\cal U}\,$ such that 
all topologies $\;\tau[{\cal F}]\;({\cal F}\in{\cal U})\;$
are homeomorphic and all
topologies $\;\sigma[{\cal F}]\;({\cal F}\in{\cal U})\;$
are homeomorphic. Indeed, put 
$\;{\cal U}\,:=\,\{\,{\cal F}_k\;|\;k\,=\,0,1,2,...\,\}\;$
where $\;{\cal F}_k\,:=\,\{\,k+S\;|\;S\in{\cal F}_0\,\}\;$
for every integer $\,k\geq 0\,$.         
Clearly, $\;{\cal F}_m\,=\,\{\,(m-n)+S\;|\;S\in{\cal F}_n\,\}\;$
whenever $\,n,m\geq 0\,$ and each family $\,{\cal F}_k\,$
is a free ultrafilter on $\,{\Bbb Z}\,$. 
We have $\,{\cal F}_n\not={\cal F}_m\,$ whenever $\,0\leq n<m\,$
because firstly precisely one of the congruence classes modulo $\,2m\,$ 
lies in $\,{\cal F}_n\,$. (Note that a union of finitely many sets 
lies in an ultrafilter only if one of these sets lies in the ultrafilter.)
And secondly, if a congruence class $\,A\,$ modulo $\,2m\,$ 
lies in $\,{\cal F}_n\,$ then the congruence class $\,(m-n)+A\,$ lies 
in $\,{\cal F}_m\,$ but not in $\,{\cal F}_n\,$. (For
$\,A\,$ and $\,(m-n)+A\,$ are disjoint.)
Finally, for each $\,k\in{\Bbb N}\,$ define an increasing bijection $\,\varphi_k\,$
from $\,{\Bbb R}\,$ onto $\,{\Bbb R}\,$ so that $\;\varphi_k(0)=0\;$
and $\;\varphi_k(n)=n+k\;$ for every $\;n\,\in\,{\Bbb Z}\setminus[-k,0]\,$.
Since $\,\varphi_k\,$ is a homeomorphism from 
the Euclidean space $\,{\Bbb R}\setminus\{0\}\,$ onto itself,
by considering the open neighborhoods of $\,0\,$
it is evident that $\,\varphi_k\,$ is a homeomorphism from 
the space $\,({\Bbb R},\tau[{\cal F}_0])\,$ onto the space $\,({\Bbb R},\tau[{\cal F}_k])\,$
and also a homeomorphism from the space $\,({\Bbb R},\sigma[{\cal F}_0])\,$ onto 
the space $\,({\Bbb R},\sigma[{\cal F}_k])\,$.
\bigskip
{\bff 6. Counting Polish spaces}
\medskip
For the proof of our second main result in Section 7 we need 
the following enumeration theorem.
\medskip
{\bf Theorem 3.} {\it There is a family $\;{\cal H}\;$
of countably infinite G$_\delta$-sets in the Euclidean space $\,{\Bbb R}\,$ 
such that the size of $\,{\cal H}\,$ is $\,c\,$ and distinct
members of $\,{\cal H}\,$ are always 
non-homeomorphic subspaces of $\,{\Bbb R}\,$.}
\medskip
{\it Proof.} We work with Cantor derivatives and is enough to consider
finite derivatives. (Note in the following that 
we regard $\,{\Bbb N}\,$ to be defined in the classical way,
i.e.~$\,0\not\in{\Bbb N}\,$.)
If $\,X\,$ is a Hausdorff space and $\,A\subset X\,$ then 
the first derivative $\,A'\,$ of 
$\,A\,$ is the set of all limit points of $\,A\,$.
Further, with $\,A^{(1)}:=A'\,$, for every $\;k\,=\,2,3,4,...\;$
the $k$-th derivative $\,A^{(k)}\,$ of $\,A\,$
is given by $\;A^{(k)}\,=\,(A^{k-1)})'\,.$
Naturally, the first derivative of any set is closed. 
Consequently, $\;A^{(m)}\supset A^{(n)}\;$ whenever $\,m\leq n\,$.
\smallskip
Now, define for each $\;n\in{\Bbb N}\;$ a compact and countably infinite 
subset $\,K_n\,$ of the interval $\,[2n,2n+1]\,$ with 
$\;\min K_n=2n\;$ and $\,\max K_n=2n+1\;$
such that $\;K_n^{(n)}=\{2n+1\}\,$.
(Simply take for $\,K_n\,$ an appropriate order-isomorphic copy
of the well-ordered set of all ordinal numbers 
$\;\alpha\leq \omega^n\,$.) Thus for $\;m,n\in{\Bbb N}\;$
the derived set $\,K_n^{(m)}\,$ contains the point $\,2n+1\,$
if and only if $\,m\leq n\,$.
Furthermore, define a discrete subset $\,E_n\,$ of 
$\;]2n+1,2n+{7\over 4}]\;$ via
$\;E_n\,:=\,\big\{\,2n+1+2^{-m}+2^{-m-k}\;\,\big|\,\;m,k\,\in\,{\Bbb N}\,\big\}\,$.
For every nonempty $\,S\subset{\Bbb N}\,$ put 
$\;\;G_S\,:=\,\bigcup_{n\in S}(K_n\cup E_n)\,$.
Since  $\,G_S\,$ is the union of the closed set 
$\;\bigcup_{n\in S}K_n\;$ and the discrete set 
$\;\bigcup_{n\in S}E_n\,$, the set $\,G_S\,$ 
is a countably infinite G$_\delta$-set in $\,{\Bbb R}\,$.
Obviously, 
$\;G_S^{(m)}\,=\,\bigcup_{n\in S}K_n^{(m)}\;$ for every $\,m\in{\Bbb N}\,$.
\smallskip
If $\,\emptyset\not=S\subset{\Bbb N}\,$ then let $\,N_S\,$ denote the 
set of all $\,x\in G_S\,$ such that no neighborhood of 
the point $\,x\,$ in the space $\,G_S\,$ is compact.
By construction, $\,x\in N_S\,$ if and only if $\;x=2n+1\;$
for some $\,n\in S\,$. Hence a moment's reflection suffices to see that
\smallskip
\centerline{$\big\{\,m\in{\Bbb N}\;\,\big|\,\;
\big(G_S^{(m)}\setminus G_S^{(m+1)}\big)\cap N_S\,\not=\,\emptyset\,\big\}
\;=\;S$}
\smallskip
for each nonempty set $\,S\subset{\Bbb N}\,$.
Thus the set $\,S\,$ can always be recovered 
from the space $\,G_S\,$  purely topologically
and hence two spaces $\,G_{S_1}\,$ and $\,G_{S_2}\,$
are never homeomorphic for distinct nonempty sets $\,S_1, S_2\subset{\Bbb N}\,$.
Thus the family $\;{\cal H}\,=\,\{\,G_S\;|\;\emptyset\not=S\subset{\Bbb N}\,\}\;$
is as desired and this concludes the proof of Theorem 3.
\bigskip
{\it Remark.} Every Polish space is homeomorphic 
to a closed subspace of the product of countably infinitely many 
copies of the real line (cf.~[3] 4.3.25).
As a consequence, every uncountable Polish space is of size $\,c\,$
and the size of a family of mutually non-homeomorphic 
Polish spaces cannot exceed $\,c\,$. Therefore, by virtue of Theorem 3,
{\it there exist precisely $\,c\,$ countably infinite Polish spaces 
up to homeomorphism.} In comparison, by [7] Proposition 2 there 
exist precisely $\,c\,$ uncountable Polish spaces up to homeomorphism.
\bigskip
{\bff 7. Completely metrizable topologies}
\medskip
{\bf Theorem 4.} {\it There exist $\,c\,$ mutually non-homeomorphic 
topologies $\,\tau\,$ on $\,{\Bbb R}\,$ 
coarser than the Euclidean topology such that
$\,({\Bbb R},\tau)\,$ is completely metrizable (and hence Polish).}
\medskip
{\it Proof.} 
Let $\;{\cal H}\;$ be a family as in Theorem 3.
Our goal is to construct for each $\,H\in{\cal H}\,$ 
a real arc $\,A_H\,$ which is a 
G$_\delta$-subset of the Euclidean space $\,{\Bbb R}^3\,$
(and hence completely metrizable) so that
$\;H\times\{0\}\times\{0\}\,\subset\,A_H\;$ 
and $\,A_H\,$ and $\,A_{H'}\,$ are never homeomorphic
for distinct $\,H,H'\in{\cal H}\,$.
\smallskip
For two points $\,P,Q\,$ in the vector space $\,{\Bbb R}^3\,$ let 
$\;[P,Q]\;$ denote the closed straight segment 
which connects the points $P$ and $Q$, 
$\;[P,Q]\,=\,\{\,\lambda P+(1\!-\!\lambda)Q\;\,|
\,\;0\leq \lambda\leq 1\,\}\,$.
Furthermore, for abbreviation, put
$\;y(n)\,:=\,2^{-n}\cos 2^{-n}\;$
and $\;z(n)\,:=\,2^{-n}\sin 2^{-n}\;$
for $\,n\in{\Bbb N}\,$.
\smallskip
For every set $\;H\,=\,\{\,a_1,a_2,a_3,...\,\}\;$ 
in the family $\,{\cal H}\,$ with $\,a_i\not=a_j\,$ for $\,i\not=j\,$
we define an injective and continuous mapping 
$\;g=g_H\;$ from $\,{\Bbb R}\,$ into $\,{\Bbb R}^3\,$ by 
\smallskip
\centerline{$\;g(t)\,=\,(t\sin t,-e^{t},0)\;$ 
for every real $\,t\leq 0\,$}
\smallskip
and so that 
$\;\;g([k,k+1])\,=\,[g(k),g(k+1)]\;\;$ for every integer $\,k\geq 0\,$ where
\smallskip
\centerline{$g(0)=(0,-1,0)$ \quad and \quad $g((1)=(0,-1,1)$\quad and}
\smallskip
\centerline{$\;g(2m)\,=\,(a_m,0,0)\;\;$ and 
$\;\;g(2m+1)\,=\,(a_m,y(m),z(m))\;\;$  for every $\,m\in{\Bbb N}\,$.}
\smallskip
The injectivity of $\,g\,$ is feasible because if 
$\,E_m\,$ is the plane through 
the three points $\;g(2m),\, g(2m\!+\!1),\,g(2m\!+\!2)\;$
then $\;E_m\,\not=\,{\Bbb R}\times{\Bbb R}\times\{0\}\;$ and
$\;E_m\cap E_n\,=\,{\Bbb R}\!\times\!\{0\}\!\times\!\{0\}\;$
whenever $\,m,n\in{\Bbb N}\,$ and $\,m\not=n\,$.
\smallskip
Let $\,H\in{\cal H}\,$ and put $\;A_H\,:=\,g_H({\Bbb R})\;$
and let $\,\overline{A_H}\,$ denote
the closure of $\,A_H\,$ in the Euclidean space $\,{\Bbb R}^3\,$.  
Trivially, $\;H\times\{0\}\times\{0\}\;$ is a G$_\delta$-set
in the space $\,{\Bbb R}^3\,$ and a subspace of $\,{\Bbb R}^3\,$
homeomorphic with $\,H\,$.
Obviously, $\;\overline{A_H}\,=\,B\!\times\!\{0\}\!\times\!\{0\}\,\cup\,A_H\;$ 
for some $\,B\subset{\Bbb R}\,$. 
Hence
$\;A_H\,=\,H\!\times\!\{0\}\!\times\!\{0\}\,\cup\,(\overline{A_H}\,\cap\,
({\Bbb R}^3\setminus {\Bbb R}\!\times\!\{0\}\!\times\!\{0\}))\;$ 
is the union of a G$_\delta$-set 
and a set which is the intersection of a closed set with an open set.
Thus $\,A_H\,$ is a G$_\delta$-set in the space $\,{\Bbb R}^3\,$
and hence the Euclidean space $\,A_H\,$ is completely metrizable.
\smallskip
A moment's reflection is sufficient to see that 
$\,H\!\times\!\{0\}\!\times\!\{0\}\,$
equals the set of all points $\,a\,$ in the space $\,A_H\,$ 
where no local basis at $\,a\,$ contains only arcwise connected sets.
Therefore, the space $\,H\,$ can essentially be recovered from 
the space $\,A_H\,$ and this finishes the proof.
\bigskip
{\it Remark.} In the previous proof one cannot replace
$\,{\cal H}\,$ with a family $\,{\cal H}'\,$ 
of mutually non-homeomorphic countably infinite and {\it closed}
subspaces of the Euclidean space $\,{\Bbb R}\,$. 
Because in view of [4] Theorem 8.1 we have
$\,|{\cal H}'|\leq \aleph_1\,$ for any such family 
$\,{\cal H}'\,$ and it 
is widely known (cf.~[3]) that $\,\aleph_1<c\,$ (i.e.~the negation of
the Continuum Hypothesis) is irrefutable.
However, by applying a theorem not proved in this paper and
with a bit greater effort concerning the notations
it is not difficult 
to modify the previous proof
starting with a family $\,{\cal H}^*\,$ 
of mutually non-homeomorphic closed subspaces of $\,{\Bbb R}\,$
such that $\,|{\cal H}^*|=c\,$ and every member of $\,{\cal H}^*\,$ 
is the union of infinitely many mutually exclusive 
intervals $\,[a,b]\,$ with $\,a<b\,$.
(Such a family $\,{\cal H}^*\,$ exists by [6] Theorem 1.)
\medskip
\bigskip
{\bff 8. Completely normal spaces of first category}
\medskip
\medskip
{\bf Theorem 5.} {\it There exist $\,2^c\,$ 
mutually non-homeomorphic 
topologies $\,\tau\in{\cal L}\,$  
such that $\,({\Bbb R},\tau)\,$ is a completely normal space of first category.}
\medskip
{\it Proof.} Let $\;B\;$ be an injective mapping from $\,{\Bbb Z}\,$ 
into the power set of $\,{\Bbb R}^3\,$ such that $\,B(k)\,$ 
is always a nonempty open ball in the Euclidean metric space $\,{\Bbb R}^3\,$ 
and that $\;\{\,B(k)\;|\;k\in{\Bbb Z}\,\}\;$ is a basis 
of the Euclidean topology of $\,{\Bbb R}^3\,$.
We define a double sequence of distinct points 
\smallskip
\centerline{$......,\,P_{-3},\,P_{-2},\,P_{-1},\,
P_{0},\,P_{1},\,P_{2},\,P_{3},\,......$}
\smallskip
in $\,{\Bbb R}^3\,$ by induction.
Start with three distinct points $\,P_{-1},P_0,P_1\,$
where $\,P_{-1}\,$ does not lie in the straight line 
through $\,P_0\,$ and $\,P_1\,$.
Suppose that for $\,n\in{\Bbb N}\,$ we have already chosen $\,2n\!+\!1\,$
distinct points $\;P_k\;$ with $\,k\in{\Bbb Z}\,$ and $\,|k|\leq n\,$. 
Then choose $\;P_{n+1}\in B(n+1)\;$ and $\;P_{-n-1}\in B(-n-1)\;$
so that 
\smallskip
(i)$\,$ three distinct points in 
$\;\{\,P_k\;|\;|k|\leq n\!+\!1\,\}\;$ never lie in one straight line,
\smallskip
(ii) four distinct points in 
$\;\{\,P_k\;|\;|k|\leq n\!+\!1\,\}\;$ never lie in one plane.
\medskip
Such a choice is always possible since neither 
finitely many straight lines nor 
finitely many planes can cover any ball $\,B(k)\,$. 
\medskip
In this way we obtain a countable, dense subset 
$\;\{\,P_k\;|\;k\in{\Bbb Z}\,\}\;$ of the Euclidean space $\,{\Bbb R}^3\,$ 
(with $\,P_k\not=P_{k'}\,$ whenever $\,k\not=k'\,$)
such that $\;[P_m,P_{m+1}]\;$ and 
$\;[P_n,P_{n+1}]\setminus\{P_n,P_{n+1}\}\;$
are disjoint whenever $\;m,n\in{\Bbb Z}\;$ and $\,m\not=n\,$.
\medskip
Now define a mapping $\,g\,$ from $\,{\Bbb R}\,$ into $\,{\Bbb R}^3\,$
so that $\;g(k)\,=\,P_k\;$ and $\,g\,$ is a continuous 
bijection from $\;[k,k+1]\;$ into $\,{\Bbb R}^3\,$ 
with $\;g([k,k+1])\,=\,[P_k,P_{k+1}]\;$ for every $\,k\in{\Bbb Z}\,$.
Then $\;g\,:\;{\Bbb R}\to{\Bbb R}^3\;$ is injective and continuous 
and hence $\,g({\Bbb R})\,$ is a real arc
within $\,{\Bbb R}^3\,$ such that $\,g({\Bbb Z})\,$ is dense in $\,{\Bbb R}^3\,$.
Therefore the Euclidean compact spaces  $\;[P_k,P_{k+1}]\;$ 
are closed subsets of the space $\,g({\Bbb R})\,$ 
whose interior in the space $\,g({\Bbb R})\,$ is empty
and hence the space $\,g({\Bbb R})\,$ is of first category.
By construction, for any nonempty open set $\,U\,$ in 
the Euclidean space $\,{\Bbb R}^3\,$ the set $\,g^{-1}(U)\,$
is an unbounded subset of $\,{\Bbb R}\,$.
Thus the topology in $\,{\cal L}\,$ corresponding with $\,g({\Bbb R})\,$
is one that satisfies the desired properties of Theorem 5.
(Moreover, the topology is  metrizable.)
\medskip
The first step is done and now we are going 
to track down $\,2^c\,$ topologies as desired.
Since $\,g({\Bbb Z})\,$ is dense in $\,{\Bbb R}^3\,$
we may fix an infinite set $\,Z\subset g({\Bbb Z})\,$ such
that $\,g(0)\in Z\,$ and
the Euclidean distance between any two points in $\,Z\,$ 
is always greater than $\,1\,$.
(In particular, $\,Z\,$ is an unbounded, countable subset of $\,{\Bbb R}^3\,$.)
Similarly as in the proof of Theorem 1,
for each of the $\,2^c\,$ free ultrafilters $\,{\cal F}\,$ on $\,Z\,$
define a topology $\,\tilde\tau[{\cal F}]\,$ on $\,{\Bbb R}^3\,$
such that $\,U\subset{\Bbb R}^3\,$ lies in the family $\,\tilde\tau[{\cal F}]\,$
if and only if $\,U\,$ is Euclidean open and satisfies $\;g(0)\not\in U\;$
or $\;U\cap Z\,\in\,{\cal F}\,$.
\medskip
Of course, by exactly the same arguments as
in the proof of Theorem 1, 
for every free ultrafilter $\,{\cal F}\,$ 
on $\,Z\,$ the topology $\,\tilde\tau[{\cal F}]\,$
is completely normal and
coarser than the Euclidean topology on $\,{\Bbb R}^3\,$
(and strictly coarser precisely at the point $\,g(0)\,$).
\medskip
Now let $\,\tau=\tilde\tau[{\cal F}]\,$ be any such topology on $\,{\Bbb R}^3\,$.
Then the set $\,g({\Bbb R})\,$ equipped with the subspace topology 
of $\,({\Bbb R}^3,\tau)\,$ is completely normal.
(Here it is essential that the property {\it completely normal} is, 
other than the property {\it normal}, hereditary.)
Since $\,g\,$ is a continuous one-to-one mapping from 
$\,({\Bbb R},\eta)\,$ into $\,({\Bbb R}^3,\tau)\,$ a fortiori,
the family $\;g^{-1}(\tau)\,:=\,\{\,g^{-1}(V)\;|\;V\in\tau\,\}\;$ 
is a topology in the family $\,{\cal L}\,$
and $\,g\,$ is a homeomorphism from 
the space $\,({\Bbb R},g^{-1}(\tau))\,$
onto the space $\,(g({\Bbb R}),\tau)\,$.
In particular, the space $\,({\Bbb R},g^{-1}(\tau))\,$ is
completely normal. Furthermore, 
every nonempty open set in the space $\,({\Bbb R},g^{-1}(\tau))\,$ 
is unbounded in $\,{\Bbb R}\,$, whence $\,({\Bbb R},g^{-1}(\tau))\,$
is a space of first category by Proposition 2.
\medskip
Trivially, $\;U\cap Z\,=\,(U\cap g({\Bbb R}))\cap Z\;$ for every Euclidean 
open set $\,U\subset{\Bbb R}^3\,$. Therefore, 
by a similar argument as in the proof of Theorem 1, 
for distinct free ultrafilters 
$\,{\cal F}_1,{\cal F}_2\,$ on $\,Z\,$ 
the relative topologies of $\,\tilde\tau[{\cal F}_1]\,$
and $\,\tilde\tau[{\cal F}_2]\,$ on the set $\,g({\Bbb R})\,$
must be distinct. (We even have 
$\,\tau_1\not\subset\tau_2\,$ for such  
distinct relative topologies $\,\tau_1,\tau_2\,$ on $\,g({\Bbb R})\,$.) 
Thus by Lemma 3 we can track down a 
family $\,{\cal U}\,$ of free ultrafilters on $\,Z\,$
such that $\,|{\cal U}|=2^c\,$ and two spaces 
$\,(g({\Bbb R}),\tilde\tau[{\cal F}_1])\,$ and $\,(g({\Bbb R}),\tilde\tau[{\cal F}_2])\,$
are never homeomorphic 
for distinct $\,{\cal F}_1,{\cal F}_2\in{\cal U}\,$.
Hence 
the topologies $\,g^{-1}(\tilde\tau[{\cal F}_1])\,$  and
$\,g^{-1}(\tilde\tau[{\cal F}_2])\,$ in the family $\,{\cal L}\,$
are never homeomorphic for distinct $\,{\cal F}_1,{\cal F}_2\in{\cal U}\,$
since $\,g\,$ is a homeomorphism from 
the space $\,({\Bbb R},g^{-1}(\tilde\tau[{\cal F}]))\,$
onto the space $\,(g({\Bbb R}),\tilde\tau[{\cal F}])\,$
for every $\,{\cal F}\in{\cal U}\,$. This concludes the proof.
\bigskip
{\bff 9. Metrizable spaces of first category}
\medskip
{\bf Theorem 6.} {\it There exist $\,c\,$ mutually non-homeomorphic 
topologies $\,\tau\in{\cal L}\,$  
such that $\,({\Bbb R},\tau)\,$ is a metrizable space of first category.}
\medskip
{\it Proof.} 
Let $\,\eta_3\,$ denote the Euclidean topology on $\,{\Bbb R}^3\,$
and for any continuous one-to-one mapping $\;g:\,{\Bbb R}\to{\Bbb R}^3\;$ let 
$\;g^{-1}(\eta_3)\,:=\,\{\,g^{-1}(V)\;|\;V\in\eta_3\,\}\;$ 
denote the topology in $\,{\cal L}\,$ corresponding with 
the real arc $\,g({\Bbb R})\,$.
Let $\,{\cal H}\,$ be a family as in Theorem 3.
Our goal is to 
construct a real arc $\,h_H({\Bbb R})\,$
within the metrizable space $\,({\Bbb R}^3,\eta_3)\,$ 
for every $\,H\in{\cal H}\,$ such that firstly $\,h_H({\Bbb Z})\,$ 
is dense in $\,{\Bbb R}^3\,$, whence every nonempty open set in the space
$\,({\Bbb R},h_H^{-1}(\eta_3))\,$ is unbounded, 
and secondly two real arcs $\,{h}_{H_1}({\Bbb R})\,$ 
and $\,h_{H_2}({\Bbb R})\,$
are never homeomorphic for distinct sets $\,H_1,H_2\in{\cal H}\,$.
\smallskip
Let $\;H\,=\,\{\,a_1,a_3,a_5,...\,\}\;$ 
be a set in the family 
$\,{\cal H}\,$ where $\;a_i\not=a_j\;$
for distinct (and always odd) indices $\,i,j\,$.
Again let $\;y(n)\,:=\,2^{-n}\cos 2^{-n}\;$
and $\;z(n)\,:=\,2^{-n}\sin 2^{-n}\;$ for $\,n\in{\Bbb N}\,$.
We firstly define $\,h=h_H\,$
on the domain $\;[0,\infty[\,$. 
Choose an injective and continuous mapping 
$\,h\,$ from $\,[0,\infty[\,$ into $\,{\Bbb R}^3\,$ 
so that 
$\;h([k,k+1])\,=\,[h(k),h(k+1)]\;$ 
for every integer $\,k\geq 0\,$ where
$\;h(k)\,=\,((-2)^{k/2},y(k),z(k))\;$ 
when $\,k\,$ is even and 
$\;h(k)\,=\,(a_{k},0,0)\;$ when $\,k\,$ is odd.
(Such a choice is clearly possible because if 
$\,E_m\,$ is the plane through 
the three points $\;h(m-1),h(m),h(m+1)\;$
for any even $\,m\geq 2\,$ 
then $\;E_m\cap E_n\,=\,{\Bbb R}\!\times\!\{0\}\!\times\!\{0\}\;$
whenever $\;2\leq m<n\,$.)
Clearly, $\,H\times\{0\}\times\{0\}\,$ is
the intersection of $\,h([0,\infty[)\,$ with the $x$-axis 
$\,{\Bbb R}\!\times\!\{0\}\!\times\!\{0\}\,$, 
and $\;h([0,\infty[)\,\cup\,{\Bbb R}\!\times\!\{0\}\!\times\!\{0\}\;$ 
is the closure of $\,h([0,\infty[)\,$ in $\,{\Bbb R}^3\,$.
\smallskip
For any Hausdorff space $\,X\,$  
let $\,W(X)\,$ denote the set of all points $\,x\,$ in
$\,X\,$ such that no local basis at $\,x\,$ 
contains only arcwise connected sets.
By construction we 
have 
\smallskip
$\;W(h([0,\infty[))\,=\,H\times\{0\}\times\{0\}\,$.
\smallskip
In view of the definition of $\,g\,$ in the proof of Theorem 5 it is plain 
to expand $\,h\,$ to  a continuous and injective mapping 
from $\,{\Bbb R}\,$ into $\,{\Bbb R}^3\,$ such that $\,h({\Bbb Z}\setminus{\Bbb N})\,$
is a dense subset of the Euclidean space $\,{\Bbb R}^3\,$.
As a consequence we have $\;W(h({\Bbb R}))\,=\,h({\Bbb R})\;$
and $\,({\Bbb R},h^{-1}(\eta_3))\,$ is a space of first category.
Moreover,
$\;W(h([t,\infty[))\,=\,H\times\{0\}\times\{0\}\;$ for every real 
$\,t\leq 0\,$ and
$\;W(h([t,\infty[))\,\subset\,H\times\{0\}\times\{0\}\;$ 
and $\;W(h(]{-\infty,t}]))\,=\,h(]{-\infty,t}])\;$
for every $\,t\in{\Bbb R}\,$.
In particular, for every $\,t\in{\Bbb R}\,$ the set
$\,W(h([t,\infty[))\,$ is countable and the set
$\,W(h(]{-\infty,t}]))\,$ is uncountable
and we have
$\;H\!\times\!\{0\}\!\times\!\{0\}\,=\,
\bigcup\,\{\;W(h([t,\infty[))\;|\;t\in{\Bbb R}\,\}\,.$
\smallskip
We finish the proof by verifying that 
$\;H\!\times\!\{0\}\!\times\!\{0\}\;$
can be recovered from the space $\,h({\Bbb R})\,$.
(Note, again, that $\;H\!\times\!\{0\}\!\times\!\{0\}\;$
and $\,H\,$ are homeomorphic.)
\smallskip
For any arcwise connected metrizable space $\,X\,$ 
let $\,{\cal Y}(X)\,$ be the family of all sets $\,Y\subset X\,$ 
such that $\,Y\,$ and $\,X\setminus Y\,$ 
are arcwise connected and 
$\,Y\setminus\{y\}\,$ is arcwise connected for some
$\,y\in Y\,$. For the Euclidean space $\,{\Bbb R}\,$ we 
clearly have $\,Y\in{\cal Y}({\Bbb R})\,$ if and only if 
$\,Y=\,]{-\infty, t}]\,$ 
or $\,Y=[t,\infty[\;$ for some $\,t\in{\Bbb R}\,$.
While for an arbitrary real arc $\,g({\Bbb R})\,$ it is not necessary 
that $\;{\cal Y}(g({\Bbb R}))\,=\,\{\,g(Y)\;|\;Y\in{\cal Y}({\Bbb R})\,\}\,$
(see the remark below), we observe that 
$\,Y\in{\cal Y}(h({\Bbb R}))\,$ 
if and only if $\;Y=h(]{-\infty, t}])\;$ 
or $\;Y=h([t,\infty[)\;$ for some $\,t\in{\Bbb R}\,$.
Therefore, $\;H\!\times\!\{0\}\!\times\!\{0\}\;$
equals the union of all sets $\,W(Y)\,$ 
where $\;Y\in{\cal Y}(h({\Bbb R}))\;$ and $\,W(Y)\,$ is countable, 
{\it q.e.d.}
\medskip
{\it Remark.} If $\,g({\Bbb R})\subset{\Bbb R}^3\,$ is a real arc
and $\,a\in{\Bbb R}\,$ 
such that $\,g(x_n)\,$ converges to $\,g(a)\,$ whenever
$\,(x_n)\,$ is an unbounded and increasing sequence of reals 
then $\;g({\Bbb R})\setminus\{g(x)\}\;$ is arcwise 
connected for every $\,x>a\,$
and $\;g([u,v])\in {\cal Y}(g({\Bbb R}))\;$
whenever $\;a<u<v\,$. 
\bigskip
{\bff 10. A complete lattice of topologies}
\medskip
As any family of topologies on a fixed set, the family $\,{\cal L}\,$
is partially ordered by the relation $\,\subset\,$.
A family $\,{\cal K}\subset{\cal L}\,$ is a {\it chain}
if and only if $\,\tau_1\subset\tau_2\,$ or $\,\tau_2\subset\tau_1\,$
whenever $\,\tau_1,\tau_2\in{\cal K}\,$.
The extreme opposite of chains of topologies are families 
of mutually incomparable topologies.
(Two topologies $\,\tau_1,\tau_2\,$ are incomparable
if and only if neither $\,\tau_1\subset\tau_2\,$ nor 
$\,\tau_2\subset\tau_1\,$.)

\medskip%\smallskip
In order to prove Theorem 1 we considered topologies in
$\,{\cal L}\,$ which are coarse at precisely one point $\,a\in{\Bbb R}\,$
(with $\,a=0\,$).
Let $\;{\cal L}_0\,:=\,\{\,\tau\in{\cal L}\;|\;C(\tau)\subset\{0\}\,\}\;$
be the family of all topologies in $\,{\cal L}\,$ which are either 
coarse precisely at the point 
$\,0\,$ or equal to the Euclidean topology $\,\eta\,$.
We have $\;|{\cal L}_0|=|{\cal L}|=2^c\;$ by the proof of Theorem 1.
Whereas, naturally, the family of all topologies on the set $\,{\Bbb R}\,$
coarser than $\,\eta\,$ is a lattice with respect to the partial ordering 
$\,\subset\,$, the partially ordered 
family $\,({\cal L},\subset)\,$ 
is not a lattice. (See the remark below.) 
However, the partially ordered 
family $\,({\cal L}_0,\subset)\,$ is a lattice.
Moreover, $\,({\cal L}_0,\subset)\,$ is a complete lattice 
(with $\,\eta\,$ as its maximum) in view of the following proposition 
which also shows that for the minimum $\,\theta\,$ of the complete
lattice $\,{\cal L}_0\,$ the space $\,({\Bbb R},\theta)\,$ 
has interesting properties. 
(Recall that a partially ordered set $\,L\,$
is a complete lattice if and only if 
every nonempty subset of $\,L\,$ has an infimum and a supremum.) 
\medskip
{\bf Proposition 6.} {\it If $\;\emptyset\not={\cal S}\subset{\cal L}_0\;$
then $\;\bigcap {\cal S}\,\in\,{\cal L}_0\,$.
If $\,{\cal K}\not=\emptyset\,$ is a chain in $\,{\cal L}_0\,$
then $\,\bigcup{\cal K}\,$ is a topology in $\,{\cal L}_0\,$,
and $\,\bigcup{\cal K}\not=\eta\,$ when $\,\eta\not\in{\cal K}\,$.
If $\;\theta\,=\,\bigcap{\cal L}_0\;$ then
the Hausdorff space $\,({\Bbb R},\theta)\,$ is compact
and any locally connected, compact                   
real arc with precisely one cut point
is homeomorphic to the space $\,({\Bbb R},\theta)\,$.}
\medskip
{\it Proof.} Let $\;\emptyset\not={\cal S}\subset{\cal L}_0\,$.
The family $\;\sigma\,:=\,\bigcap{\cal S}\;$ is 
a topology on $\,{\Bbb R}\,$ coarser than $\,\eta\,$ 
since, generally, the lattice of all topologies on any set
is closed under arbitrary intersections.
The topology $\,\sigma\,$ is Hausdorff because
$\,\sigma\,$ and $\,\eta\,$ coincide on $\,{\Bbb R}\setminus\{0\}\,$
and if, say, $\,x>0\,$ then $\,0\,$ and $\,x\,$ can be separated
by the $\sigma$-open sets $\;{\Bbb R}\setminus[{x\over 3},3x]\;$ 
and $\;]{x\over 2},2x[\,$. (Since $\;[{x\over 3},3x]\;$
is $\tau$-compact for every $\,\tau\in{\cal L}\,$, the set 
$\,{\Bbb R}\setminus[{x\over 3},3x]\,$ is $\tau$-open  
for every $\,\tau\in{\cal S}\,$.)
If $\,{\cal S}\not=\{\eta\}\,$ then 
$\,C(\sigma)=\{0\}\,$ by Proposition 1. 
Hence, $\,\sigma\in{\cal L}_0\,$.
Recall that 
if $\,\tau\in{\cal L}_0\,$ and $\,0\in U\in\tau\,$ and $\,V\in\eta\,$ 
then $\,U\cup V\in\tau\,$. And, by Proposition 1, $\;]{-1,1}[\,\in\tau\;$ for
$\,\tau\in{\cal L}_0\,$ only if $\,\tau=\eta\,$.
Consequently, 
the family $\,\bigcup{\cal S}\,$ is closed under arbitrary unions
and we have $\,\bigcup{\cal S}\not=\eta\,$ when $\,\eta\not\in{\cal S}\,$.
And if $\,{\cal S}\,$ is a chain then $\,\bigcup{\cal S}\,$ 
is closed under finite intersections and hence 
$\,\bigcup{\cal S}\,$ is a topology on $\,{\Bbb R}\,$ coarser than $\,\eta\,$ 
and finer than the Hausdorff topology $\,\bigcap{\cal S}\,$,
whence $\,\bigcup{\cal S}\in{\cal L}_0\,$.
\smallskip
Define a topology $\,\tau_0\in{\cal L}\,$ 
by declaring a set $\;U\subset{\Bbb R}\;$ 
$\tau_0$-open if and only if the set $\,U\,$ is $\eta$-open 
and either $\,0\not\in U\,$
or $\;U\,\supset\,\{0\}\cup({\Bbb R}\setminus[-t,t])\;$ for some $\,t>0\,$.
Then $\,C(\tau_0)=\{0\}\,$ and hence $\,\tau_0\in{\cal L}_0\,$.
Let $\,K\,$ be the union of two congruent circles 
in the plane $\,{\Bbb R}^2\,$ which meet in precisely one point.
Then $\,K\,$ (which looks like the digit $8$ 
or the symbol $\infty$) is 
an arcwise connected and locally arcwise connected compact subspace 
of the Euclidean plane $\,{\Bbb R}^2\,$ with precisely one cut point.
(Recall that $\,x\,$ is a cut point of a connected space
$\,X\,$ if and only if $\,X\setminus\{x\}\,$ is not connected.) 
It is immediately obvious that $\,K\,$ is a real arc 
which is homeomorphic to the space $\,({\Bbb R},\tau_0)\,$.
(Of course, $\,0\,$ is the unique cut point in the 
arcwise connected space $\,({\Bbb R},\tau_0)\,$.) 
It is well-known that
any locally connected, compact 
real arc with precisely one cut point is homeomorphic to $\,K\,$ (cf.~[9]).
Finally, the topologies $\,\tau_0\,$ and $\,\bigcap{\cal L}_0\,$
must be identical because $\,\tau_0\in{\cal L}_0\,$
and $\,\tau_0\subset\tau\,$ for every $\,\tau\in{\cal L}_0\,$
since if $\;0\in U\in\tau_0\;$ then $\;{\Bbb R}\setminus U\;$ is
Euclidean compact and hence $\tau$-closed for
every $\,\tau\in{\cal L}_0\,$, {\it q.e.d.}
\medskip
{\it Remark.} If $\,a\in{\Bbb R}\,$ and $\,\varphi_a(x)=x+a\,$ for every 
$\,x\in{\Bbb R}\,$ and
$\,\tau_0\in{\cal L}_0\,$ is compact then
$\;\tau_a\,:=\,\{\,\varphi_a(U)\;|\;U\in\tau_0\,\}\;$ is a 
topology in $\,{\cal L}\,$ with $\,C(\tau_a)=\{a\}\,$ and 
hence $\,\tau_a\not=\tau_{a'}\,$ whenever $\,a\not=a'\,$.
Each topology $\,\tau_a\,$ 
is compact since $\,\varphi_a\,$ is a homeomorphism from 
$\,({\Bbb R},\tau_0)\,$ onto $\,({\Bbb R},\tau_a)\,$. Thus by Proposition 6,
$\,{\cal L}\,$ contains $\,c\,$ (homeomorphic) compact topologies.
Therefore, the partially ordered family $\,({\cal L},\subset)\,$
is not a lattice because
if $\,\tau,\tau'\,$ are distinct compact topologies in $\,{\cal L}\,$
then $\,\{\tau,\tau'\}\,$
has no infimum in $\,({\cal L},\subset)\,$
since a topology cannot be T$_2$ if it is 
strictly coarser than a T$_2$-compact topology.
(In particular, every nonempty chain of compact topologies
in $\,{\cal L}\,$ is a singleton.)
It is also worth mentioning that
if for $\,\tau\in{\cal L}\,$ the space $\,({\Bbb R},\tau)\,$ is
compact then it must be second countable. Because, naturally, the sets
$\;]r_1,r_2[\;$ with $\;r_1,r_2\in{\Bbb Q}\;$ form a network of 
$\,\tau\,$ and (cf.~[2] 3.3.5.) 
any compact Hausdorff space has a countable basis
if it has a countable network.

\bigskip
{\bff 11. Long chains of homeomorphic topologies}
\medskip
The topologies in the family $\,{\cal T}\subset{\cal L}\,$ constructed in the 
proof of Theorem 1  are mutually non-homeomorphic
and mutually incomparable.
If $\,\tau_z\in{\cal L}\,$ are the completely metrizable topologies 
defined by the real arcs $\,g_z({\Bbb R})\,$
in the proof of Proposition 5 
then $\;\{\,\tau_z\;|\;z\in{\Bbb R}\,\}\;$ is a family of homeomorphic and
mutually incomparable topologies.
(They are mutually incomparable because if $\,r,s\in{\Bbb R}\,$ and $\,r\not=s\,$
then the sequence $\,(1+r+\pi n)\,$ converges to $\,r\,$ in the space 
$\,({\Bbb R},\tau_r)\,$, whereas in the space $\,({\Bbb R},\tau_s)\,$
the same sequence converges to $\,s\,$ when
$\,{r-s\over\pi}\in{\Bbb Z}\,$ and diverges
when $\,{r-s\over\pi}\not\in{\Bbb Z}\,$.)
However, a simple modification of the real arc $\,g_z({\Bbb R})\,$ 
makes it possible to track down a chain of homeomorphic topologies in 
$\,{\cal L}\,$.                                                           
\medskip
{\bf Proposition 7.} {\it There exists a chain 
$\,{\cal J}\subset{\cal L}\,$
such that $\,|{\cal J}|=c\,$ and all spaces $\,({\Bbb R},\tau)\,$
with $\,\tau\in{\cal J}\,$
are completely metrizable and homeomorphic.}
\smallskip
{\it Proof.} For $\,z\in{\Bbb R}\,$ consider the mapping 
$\;g_z:\,{\Bbb R}\to{\Bbb R}^2\;$ from the proof of Proposition 5 and for
$\,-1<a<0\,$ put $\;\tilde g_a(t)=g_0(t)\;$ when $\,t\geq 0\,$ 
and $\;\tilde g_a(t)=(0,-t)\;$ when $\;a\leq t\leq 0\;$
and $\;\tilde g_a(t)=(t-a,-a)\;$ when $\;t\leq a\,$.
For $\,-1<a<0\,$ 
let $\,\tilde\tau_a\,$ be the topology in $\,{\cal L}\,$ 
corresponding with the Euclidean continuous injective mapping 
$\;\tilde g_a:\,{\Bbb R}\to{\Bbb R}^2\,$. 
Then $\;C(\tilde\tau_a)=[a,0]\;$ and
$\,({\Bbb R},\tilde\tau_a)\,$ is completely metrizable
since $\,\tilde g_a({\Bbb R})\,$ is a G$_\delta$-subset of $\,{\Bbb R}^2\,$.
Obviously, $\,\tilde\tau_r\,$ is a proper subset of $\,\tilde\tau_s\,$ 
whenever $\;-1<r<s<0\,$.
All spaces $\,({\Bbb R},\tilde\tau_a)\,$ with $\,-1<a<0\,$ are homeomorphic
because a moment's reflection suffices to see
that if $\;-1<r<s<0\;$ then there is a homeomorphism from 
the Euclidean plane $\,{\Bbb R}^2\,$ onto itself which 
maps $\,\tilde g_r({\Bbb R})\,$ onto $\,\tilde g_s({\Bbb R})\,$, {\it q.e.d.}
\medskip\smallskip
The chain $\,{\cal J}\,$ of homeomorphic topologies 
constructed in the previous proof
is disjoint from the lattice $\,{\cal L}_0\,$.
If $\,{\cal T}\,$ is a family as in Theorem 1 then 
$\,{\cal T}\subset{\cal L}_0\,$ 
but there is no chain $\,{\cal K}\subset{\cal T}\,$ with $\,|{\cal K}|>1\,$.
Nevertheless, the following theorem shows that the lattice $\,{\cal L}_0\,$
contains very long chains of homeomorphic topologies.
(In the following, as usual,
if $\,{\cal K}_2\,$ is a $\subset$-chain
and $\,{\cal K}_1\subset{\cal K}_2\,$ then 
$\,{\cal K}_1\,$ is dense in $\,{\cal K}_2\,$
if and only if for every pair $\;X,Y\in{\cal K}_2\;$ 
with $\,X\subset Y\,$ and $\,X\not=Y\,$ there exists 
a set $\,Z\,$ in $\;{\cal K}_1\setminus\{X,Y\}\;$
such that $\;X\subset Z\subset Y\,$.)
\medskip
{\bf Theorem 7.} {\it  The lattice $\,{\cal L}_0\,$
contains four chains 
$\,{\cal K}_0,{\cal K}_1,{\cal K}_2,{\cal K}_3\,$
of (the maximal possible) size $\,c\,$ such that for $\,i\in\{0,1,2,3\}\,$
all spaces $\,({\Bbb R},\tau)\,$ with $\,\tau\in{\cal K}_i\,$
are homeomorphic, and
\smallskip
{\rm (i)} if $\,\tau\in{\cal K}_0\,$ then the space $\,(R,\tau)\,$ 
is second countable but not regular,
\smallskip
{\rm (ii)} if $\,\tau\in{\cal K}_1\,$ then the space $\,(R,\tau)\,$ 
is neither regular nor first countable, 
\smallskip
{\rm (iii)} if $\,\tau\in{\cal K}_2\,$ then the space $\,(R,\tau)\,$ 
is completely normal but not first countable,
\smallskip
{\rm (iv)} if $\,\tau\in{\cal K}_3\,$ then the space $\,(R,\tau)\,$ 
is completely metrizable,
\smallskip
{\rm (v)} ${\cal K}_0\cup{\cal K}_1\cup{\cal K}_2\;$
is a chain and $\,{\cal K}_i\,$ is 
dense in $\,{\cal K}_0\cup{\cal K}_1\cup{\cal K}_2\,$ for every 
$\,i\in\{0,1,2\}\,$,
\smallskip
{\rm (vi)} every topology in $\,{\cal K}_0\cup{\cal K}_1\cup{\cal K}_2\,$
is coarser than every topology in $\,{\cal K}_3\,$.}
\medskip
{\it Proof.} The size of $\,{\cal K}_i\,$
cannot exceed $\,c\,$ by Lemma 2.
In order to obtain a chain $\,{\cal K}_3\,$ as desired,
for real $\,\alpha\geq 0\,$ 
define an injective and Euclidean continuous 
mapping $\,h_\alpha\,$ from $\,{\Bbb R}\,$ into $\,{\Bbb R}^2\,$
by 
$\;h_\alpha(t)=(t,-t)\;$ for $\,t \leq 1\,$ and 
$\;h_\alpha(t)=(1,t-2)\;$ for $\,1\leq t\leq 2\,$ and 
$\;h_\alpha(t)=(2t^{-1},t^\alpha |\sin(\pi t)|)\;$ 
for $\,t\geq 2\,$.                       
\smallskip
Obviously $\,h_\alpha({\Bbb R})\,$ is a G$_\delta$-subset of $\,{\Bbb R}^2\,$
for every $\,\alpha\geq 0\,$.
All sets $\,h_\alpha({\Bbb R})\,$ with $\,\alpha\geq 0\,$
are homeomorphic subspaces of $\,{\Bbb R}^2\,$ because for every $\,\alpha\geq 0\,$
the mapping  $\;(t,h_0(t))\mapsto(t,h_\alpha(t))\;$ 
with $\,t\,$ running through $\,{\Bbb R}\,$ is clearly 
a homeomorphism from the real arc $\,h_0({\Bbb R})\,$
onto the real arc $\,h_\alpha({\Bbb R})\,$.
Let $\,\mu[\alpha]\,$ be the topology in $\,{\cal L}\,$ 
corresponding with $\,h_\alpha\,$. Thus
$\,\mu[\alpha]\in{\cal L}_0\,$  and in the space $\,({\Bbb R},\mu[\alpha])\,$ the family 
$\;\{\,B(\alpha,\varepsilon)\;|\;\varepsilon>0\,\}\;$ is 
a local basis at the point $\,0\,$ where

\centerline{$\;B(\alpha,\varepsilon)\;:=\;\,]{-\varepsilon,\varepsilon}[\;\cup\,
\;\{\,t\in{\Bbb R}\;|\;t>{2\over \varepsilon}\;\land\;
t^{\alpha}|\sin(\pi t)|<\varepsilon\,\}\,.$}
\smallskip
(Obviously, 
$\;h_{\alpha}^{-1}(]{-\varepsilon,\varepsilon}[^2\cap h_{\alpha}({\Bbb R}))\,= 
B(\alpha,\varepsilon)\;$ for every positive $\,\varepsilon<1\,$.)
If $\;0\leq\alpha_1\leq \alpha_2\;$ then 
$\;B(\alpha_1,\varepsilon)\supset B(\alpha_2,\varepsilon)\;$
for every $\,\varepsilon>0\,$ 
and hence $\,\mu[\alpha_1]\subset\mu[\alpha_2]\,$.
If $\;0\leq\alpha_1<\alpha_2\;$ then 
$\,\mu[\alpha_1]\not=\mu[\alpha_2]\,$
because the $\mu[\alpha_2]$-open set
$\,B(\alpha_2,1)\,$ 
cannot be $\mu[\alpha_1]$-open since it is plain 
that 
$\;B(\alpha_1,\varepsilon)\not\subset B(\alpha_2,1)\;$ 
for every $\,\varepsilon>0\,$.
So we define $\;{\cal K}_3\,:=\,\{\,\mu[\alpha]\;|\;\alpha\geq 0\,\}\,$.
\medskip
In order to find appropriate chains $\;{\cal K}_0,{\cal K}_1,{\cal K}_2\;$
we define a family $\,{\cal D}\subset{\cal L}_0\,$ 
so that the partially ordered set $\;({\cal D},\subset)\;$
is a Boolean algebra isomorphic with the power set of $\,{\Bbb R}\,$.
Write $\;x+Y\,:=\,\{\,x+y\;|\;y\in Y\,\}\;$ 
for $\,x\in{\Bbb R}\,$ and $\,Y\subset{\Bbb R}\,$.
For any set $\;D\subset[-{1\over 2},{1\over 2}[\;$ define a topology 
$\,\tau(D)\in{\cal L}\,$ by 
declaring $\,U\subset{\Bbb R}\,$ open if and only if 
$\,U\,$ is Euclidean open and either $\;0\not\in U\;$ or
$\;U\,\supset\,\{0\}\cup\bigcup_{k=n}^\infty k+D\;$ for some 
$\,n\in{\Bbb N}\,$. It is plain that $\,\tau(D)\,$
is a well-defined topology on $\,{\Bbb R}\,$ and that 
$\,\tau(D)\in{\cal L}_0\,$. 
\smallskip
Obviously,
$\,\tau(\emptyset)=\eta\,$ and 
$\;\tau(B)\subset\tau(A)\;$ whenever $\;A\subset B\subset[-{1\over 2},{1\over 2}[\,$.
Furthermore $\;\tau(A)\not=\tau(B)\;$ 
when $\,A,B\,$ are distinct subsets of $\,[-{1\over 2},{1\over 2}[\,$.
Moreover, if $\,B\not\subset A\,$ then 
$\,\tau(A)\not\subset\tau(B)\,$. 
Because if $\;z\in B\setminus A\;$ then it is clear that 
the Euclidean open
set $\;{\Bbb R}\setminus(z+{\Bbb N})\;$ lies in $\,\tau(A)\,$ but not in $\,\tau(B)\,$.
Therefore, if 
\smallskip
\centerline{$\;{\cal D}\,:=\,
\{\,\tau(D)\;|\;D\subset[-{1\over 2},{1\over 2}[\,\}\;$}
\smallskip
and $\,g\,$ is a bijection from $\,{\Bbb R}\,$ onto $\,[-{1\over 2},{1\over 2}[\;$
then $\;X\,\mapsto\,\tau([-{1\over 2},{1\over 2}[\setminus g(X))\;$
is an isomorphism 
from the Boolean algebra of all subsets of $\,{\Bbb R}\,$
onto the partially ordered set $\,({\cal D},\subset)\,$.
\smallskip
A moment's reflection suffices to see that $\,\tau(D)\subset\mu[\alpha]\,$
for every $\,\alpha\geq 0\,$ if $\,D\subset[-{1\over 2},{1\over 2}[\;$
and $\,0\,$ is an interior point of $\,D\,$ in the Euclidean space $\,{\Bbb R}\,$.
Therefore, in order to achieve (vi) we choose mutually disjoint sets 
$\,\Lambda_0,\Lambda_1,\Lambda_2\subset\,]0,{1\over 3}[\;$ 
of size $\,c\,$ which are dense in 
$\;]0,{1\over 3}[\;$ and define 
$\;{\cal K}_0\,:=\,\{\,\tau([-\lambda,\lambda])\;|\; 
\lambda\in\Lambda_0\,\}\;$ and
$\;{\cal K}_1\,:=\,\{\,\tau([-\lambda,\lambda[)\;|\;
\lambda\in\Lambda_1\,\}\;$ and
$\;{\cal K}_2\,:=\,\{\,\tau(]{-\lambda,\lambda}[)\;|\;
\lambda\in \Lambda_2\,\}\,$.
The specific choice of 
$\,\Lambda_0,\Lambda_1,\Lambda_2\,$
is made for saving the 
density condition (v) because 
if $\;A\subset B\subset[-{1\over 2},{1\over 2}[\;$
and $\,|B\setminus A|=1\,$ then no topology from $\,{\cal D}\,$
lies strictly between $\,\tau(B)\,$ and $\,\tau(A)\,$.
Clearly, if $\;0<\lambda,\lambda'<{1\over 3}\;$
and $\,f\,$ is any strictly increasing function 
from $\,{\Bbb R}\,$ onto $\,{\Bbb R}\,$ 
with $\,f(0)=0\,$ and 
$\;f(n\pm \lambda)\,=\,n\pm \lambda'\;$ for every $\,n\in{\Bbb N}\,$ 
then $\,f\,$ is a homeomorphism from 
$\,({\Bbb R},\tau([-\lambda,\lambda]))\,$  onto $\,({\Bbb R},\tau([-\lambda',\lambda']))\,$
and from 
$\,({\Bbb R},\tau([-\lambda,\lambda[))\,$  onto $\,({\Bbb R},\tau([-\lambda',\lambda'[))\,$
and from 
$\,({\Bbb R},\tau(]{-\lambda,\lambda}[))\,$  
onto $\,({\Bbb R},\tau(]{-\lambda',\lambda'}[))\,$.
So the definitions of the four chains $\,{\cal K}_i\,$
do the job provided that (i) and (ii) and (iii) hold.
\medskip
For $\,T\subset{\Bbb R}\,$ put
$\;\Gamma(T)\,:=\,\{\,e^{2\pi it}\;|\;t\in T\,\}\,$.
So $\;\Gamma({\Bbb R})=\Gamma([-{1\over 2},{1\over 2}[)\;$ is the unit circle
$\;x^2+y^2=1\;$ in $\,{\Bbb R}^2\,$ and $\,\Gamma(D)\subset\Gamma({\Bbb R})\,$
for $\,D\subset[-{1\over 2},{1\over 2}[\,$. 
We finish the proof by verifying the 
nice observation that
for every $\,D\subset[-{1\over 2},{1\over 2}[\,$,
\smallskip
(1)\quad $\,({\Bbb R},\tau(D))\,$ {\it is second  countable if and only if 
$\,\Gamma(D)\,$ is open in $\,\Gamma({\Bbb R})\,$},
\smallskip
(2)\quad $\,({\Bbb R},\tau(D))\,$ {\it is regular if and only if 
$\,\Gamma(D)\,$ is closed in $\,\Gamma({\Bbb R})\,$},
\smallskip
Note that by Lemma 4 the space
$\,({\Bbb R},\tau(D))\,$ is regular if and only if $\,({\Bbb R},\tau(D))\,$ 
is completely normal.
\smallskip
If $\,\Gamma(D)\,$ is open in $\,\Gamma({\Bbb R})\,$ then 
$\;\{\,]{-n^{-1},n^{-1}}[\,\cup\bigcup_{k=n}^\infty k+D\;|
\;n\in{\Bbb N}\,\}\;$
is clearly a local basis at $\,0\,$ in the space 
$\,({\Bbb R},\tau(D))\,$, whence $\,({\Bbb R},\tau(D))\,$ is second countable 
by Lemma 4.
If $\,\Gamma(D)\,$ is not closed in $\,\Gamma({\Bbb R})\,$ then 
for some $\;b\,\in\,[-{1\over 2},{1\over 2}[\setminus D\;$ 
the point $\,e^{2\pi ib}\,$ is a limit point of $\,\Gamma(D)\,$
in $\,\Gamma({\Bbb R})\,$.
So the Euclidean closed set $\;b+{\Bbb N}\;$ is $\tau(D)$-closed
and, obviously, the point $\,0\,$ and the set $\;b+{\Bbb N}\;$
can not be separated by $\tau(D)$-open sets, whence $\,\tau(D)\,$
is not regular. If $\,\Gamma(D)\,$ is closed in $\,\Gamma({\Bbb R})\,$ then, 
by the same arguments as in the proof of Theorem 1,
the space $\,({\Bbb R},\tau(D))\,$  
is regular. 
(One can adopt the proof line by line
with the only modification that the set $\,B=\{0\}\cup({\Bbb Z}\setminus A)\,$
is replaced by $\,B=\{0\}\cup\bigcup_{n=k}^\infty n+D\,$
where $\,k\in{\Bbb N}\,$ is chosen so that 
$\,A\cap (n+D)=\emptyset\,$ whenever $\,n\geq k\,$.)
\smallskip
Finally, assume that $\,\Gamma(D)\,$ is not open in $\,\Gamma({\Bbb R})\,$
and choose $\,d\in D\,$ so that $\,e^{2\pi i d}\,$ is not an interior point
of $\,\Gamma(D)\,$ in $\,\Gamma({\Bbb R})\,$.
Suppose that a countable family $\;\{\,B_1,B_2,B_3,...\,\}\;$
of Euclidean open sets 
is a local basis at $\,0\,$ in the space $\,({\Bbb R},\tau(D))\,$.
Let $\,k_1\,$ be the least positive integer $\,n\,$ such that 
$\,B_1\supset n+D\,$. If $\,k_m\,$ is already 
defined then let $\,k_{m+1}\,$ be the least integer 
$\;n>k_{m}\,$ such that $\;B_{m+1}\supset n+D\,$.
For every $\,m\in{\Bbb N}\,$ choose a small $\,\epsilon_m>0\,$  
such that $\;\Gamma(]d-\epsilon_m,d+\epsilon_m[)\not\subset \Gamma(D)\;$
and $\;]k_m+d-\epsilon_m,k_m+d+\epsilon_m[\,\subset B_{m}\,$.
Then for every $\,m\in{\Bbb N}\,$ we can choose a point 
$\,x_m\,$ in $\;]k_m+d-\epsilon_m,k_m+d+\epsilon_m[\,\setminus(k_m+D)\,.$
Then the set $\;V\,:=\,{\Bbb R}\setminus\{\,x_m\;|\;m\in{\Bbb N}\,\}\;$ is $\tau(D)$-open 
and hence $\,V\supset B_n\,$ for some $\,n\in{\Bbb N}\,$.
So we obtain the contradiction that 
$\;x_n\in B_n\subset V\;$ and $\;x_n\not\in V\,$ for some $\,n\in{\Bbb N}\,$.
Thus the assumption on $\;\{\,B_1,B_2,B_3,...\,\}\;$ is false
and hence $\,\tau(D)\,$ is not first countable. 
This concludes the proof of Theorem 7.
\medskip
{\it Remark.} The maximum of
the Boolean algebra $\,({\cal D},\subset)\,$
is $\,\tau(\emptyset)=\eta\,$.
The topology $\;\tau([-{1\over 2},{1\over 2}[)\;$ is the minimum 
of $\,{\cal D}\,$ and it is plain that 
$\,({\Bbb R},\tau([-{1\over 2},{1\over 2}[))\,$ is homeomorphic to 
the subspace $\;\Gamma^*\,:=\,\Gamma({\Bbb R})\,\cup\,\{0\}\times[1,\infty[\;$
of the Euclidean plane $\,{\Bbb R}^2\,$.
It is well-known 
that any locally connected, locally compact but not compact
real arc is homeomorphic either to $\,\Gamma^*\,$ 
or to the real line (cf.~[9]).
In view of (1) and (2),  
the maximum and the minimum of the Boolean algebra $\,{\cal D}\,$ are 
the only metrizable topologies in $\,{\cal D}\,$.
In view of (2) and \hbox{Lemma 4} and $\,|\eta|=c\,$, 
precisely $\,c\,$ topologies in 
$\,{\cal D}\,$ are completely normal, whence the proof of Theorem 1 
is not dispensable. On the contrary, in view of Lemma 3 and Proposition 4 
and the well-known fact that $\,{\Bbb R}^2\,$ has only $\,c\,$ 
Euclidean closed subsets (and the trivial fact that $\,\Gamma({\Bbb R})\,$
has $\,2^c\,$ subsets), an alternative proof of Theorem 2
(which does not use ultrafilters) is provided by (2).
\bigskip
{\bff 12. Countably generated topologies}
\medskip\smallskip
Only $\,c\,$ topologies in the Boolean algebra $\,{\cal D}\,$
are first countable.
But all topologies in $\,{\cal D}\,$ 
satisfy an interesting  countability condition
weaker than first countability.
Let $\,{\cal L}_0^*\,$ denote the family of all topologies in $\,{\cal L}_0\,$
such that  $\;{\cal N}_\tau(0)\,=\,{\cal N}_\eta(0)\cap{\cal F}\;$ 
for some filter $\,{\cal F}\,$ on $\,{\Bbb R}\,$ which is generated by a countable 
filter base. In other words, there is a countable filter base $\,{\cal B}\,$
of subsets of $\,{\Bbb R}\,$ such that 
$\;\eta\cap{\cal N}_\tau(0)\,=\,
\{\,U\in\eta\;|\;\exists\,B\in{\cal B}:\;U\,\supset\,\{0\}\cup B\,\}\,$.
So if $\,\tau\in{\cal L}_0\,$ is first countable then $\,\tau\in{\cal L}_0^*\,$.
The converse is not true since
$\,{\cal D}\subset{\cal L}_0^*\,$. In particular, $\;|{\cal L}_0^*|=|{\cal D}|=2^c\,$.
Whereas for $\;A\subset B\subset[-{1\over 2},{1\over 2}[\;$
with $\,|B\setminus A|=1\,$ 
there is no topology $\,\tau\in{\cal D}\,$ 
strictly between $\,\tau(B)\,$ and $\,\tau(A)\,$,
the following theorem implies that 
between $\,\tau(B)\,$ and $\,\tau(A)\,$ there 
lie $\,c\,$ comparable and $\,c\,$ incomparable 
topologies from $\,{\cal L}_0^*\,$
and also $\,2^c\,$
incomparable topologies from $\;{\cal L}_0\setminus{\cal L}_0^*\,$.
\medskip       
{\bf Theorem 8.} {\it If $\,\tau_1\in{\cal L}_0^*\,$ 
is strictly coarser than $\,\tau_2\in{\cal L}_0\,$ 
then there are a chain $\,{\cal R}\subset{\cal L}_0\,$ 
with $\,|{\cal R}|=c\,$ and two families  
$\,{\cal S}\subset {\cal L}_0\,$ and 
$\;{\cal T}\,\subset\,{\cal L}_0\setminus{\cal L}_0^*\,$ 
of mutually incomparable topologies 
with $\,|{\cal S}|=c\,$ and $\,|{\cal T}|=2^c\,$ 
such that $\;\tau_1\subset\tau\subset\tau_2\;$
for every $\;\tau\,\in\,{\cal R}\cup{\cal S}\cup{\cal T}\,$.
Additionally $\;{\cal R},{\cal S}\,\subset\,{\cal L}_0^*\;$
can be achieved if $\,\tau_2\in{\cal L}_0^*\,$.
For $\,\tau_2=\eta\,$ it can be achieved that        
$\;{\cal R},{\cal S}\,\subset\,{\cal L}_0^*\;$ and all topologies in 
$\,{\cal R}\cup{\cal S}\,$ are homeomorphic.}
\medskip
{\it Proof.} First of all, if $\;\eta\cap{\cal N}_{\tau}(0)\,=\,
\{\,U\in\eta\;|\;\exists\,B\in{\cal B}:\;U\,\supset\,\{0\}\cup B\,\}\;$
for $\,\tau\in{\cal L}_0\,$ and a filter base $\,{\cal B}\,$
then 
$\;\eta\cap{\cal N}_{\tau}(0)\,=\,
\{\,U\in\eta\;|\;\exists\,B\in{\cal B}:\;U\,\supset\,
\{0\}\cup(B\setminus[{-1,1}])\,\}\,$.
Indeed, if $\,U\in\eta\,$ contains    
$\;\{0\}\cup(B_1\setminus[{-1,1}])\;$ for some $\,B_1\in{\cal B}\,$ 
then $\,U\,$ contains $\;]{-k^{-1},k^{-1}}[\,\cup(B_1\setminus[{-1,1}])\;$
for some $\,k>1\,$.
Since $\;V_k\,:=\,{\Bbb R}\setminus([{-k,-k^{-1}}]\cup[k^{-1},k])\;$ lies in 
$\,\eta\cap{\cal N}_{\tau}(0)\,$, we have $\,B_2\subset V_k\,$
for some $\,B_2\in{\cal B}\,$ 
and hence $\;U\,\supset\, B_1\cap B_2\,$.
Thus, since $\,{\cal B}\,$ is a filter base, 
we have $\;B\,\subset\, B_1\cap B_2\,\subset\,U\;$
for some $\,B\in{\cal B}\,$.
There is an important consequence of the 
two representations of $\;\eta\cap{\cal N}_{\tau}(0)\,$.
If $\,\eta\not=\tau\in{\cal L}_0\,$
and a filter base $\,{\cal B}\,$ generates a filter $\,{\cal F}\,$
with $\;{\cal N}_\tau(0)\,=\,{\cal N}_\eta(0)\cap{\cal F}\;$ 
then the family 
$\;{\cal B}'\,:=\,\{\,B\setminus[-1,1]\;|\;B\in{\cal B}\,\}\;$ 
does not contain $\,\emptyset\,$ 
and hence $\,{\cal B}'\,$ is a filter base 
which generates a filter $\,{\cal F}'\,$
with $\;{\cal N}_\tau(0)\,=\,{\cal N}_\eta(0)\cap{\cal F}'\,$. 
\smallskip
Let $\,\tau_1\in{\cal L}_0^*\,$ be
a proper subset of $\,\tau_2\in{\cal L}_0\,$.
Let $\,{\cal B}_1\,$ and $\,{\cal B}_2\,$ be families 
of subsets of $\;{\Bbb R}\setminus[-1,1]\;$ 
such that $\,{\cal B}_1\,$ is a countable filter base 
and $\,{\cal B}_2\,$ is a filter base when $\,\tau_2\not=\eta\,$
and $\,{\cal B}_2=\{\emptyset\}\,$ when $\,\tau_2=\eta\,$
and $\;\eta\cap{\cal N}_{\tau_i}(0)\,=\,
\{\,U\in\eta\;|\;\exists\,B\in{\cal B}_i:\;U\,\supset\,\{0\}\cup B\,\}\;$
for $\,i\in\{1,2\}\,$.
We may assume that $\;{\cal B}_1=\{\,A_1,A_2,A_3,...\,\}\;$ 
where $\,A_n\,$ is a proper subset of $\,A_m\,$ 
whenever $\,m<n\,$. Since $\,\tau_1\,$ is strictly coarser
than $\,\tau_2\,$, we can fix $\,D\in{\cal B}_2\,$
such that $\,A_n\not\subset D\,$ for every $\,n\in{\Bbb N}\,$.
Since for every $\,k\in{\Bbb N}\,$ we have
$\;A_n\subset V_k\;$ and hence $\;A_n\subset {\Bbb R}\setminus[-k,k]\;$
for some $\,n\in{\Bbb N}\,$,
we can choose a sequence $\;a_1,a_2,a_3,...\;$ of distinct reals 
such that always $\;a_n\,\in\,A_n\setminus D\;$ and either $\;a_n>n\;$
for every $\,n\in{\Bbb N}\,$ or $\;a_n<-n\;$ for every $\,n\in{\Bbb N}\,$.
Then $\;\{\,a_1,a_2,a_3,...\,\}\;$ is
 disjoint from $\;D\cup[-1,1]\;$ and Euclidean closed and discrete.
Consequently, every subset of $\;\{\,a_1,a_2,a_3,...\,\}\;$ is
$\tau_2$-closed.           
\smallskip                                  
For every infinite set $\,S\subset{\Bbb N}\,$ 
define a topology $\,\rho[S]\in{\cal L}_0\,$ with $\,\rho[S]\subset\tau_2\,$
so that an $\tau_2$-open neighborhood $\,U\,$ of $\,0\,$ 
is $\rho[S]$-open if and only if 
$\;U\,\supset\,\{\,a_n\;|\;k\leq n\in S\,\}\;$
for some  $\,k\in{\Bbb N}\,$. We have $\,\tau_1\subset\rho[S]\,$
since $\;\{\,a_n\;|\;n\geq k\,\}\subset A_k\;$ for every $\,k\in{\Bbb N}\,$.
Obviously, $\;\rho[S_1]\subset\rho[S_2]\;$ when $\,S_1\supset S_2\,$.
Furthermore, if $\,S_2\setminus S_1\,$ is an infinite set
then $\;\rho[S_1]\not\subset\rho[S_2]\;$ 
because the $\tau_2$-open set $\;{\Bbb R}\setminus\{\,a_n\;|\;n\not\in S_1\,\}\;$
is $\rho[S_1]$-open but not $\rho[S_2]$-open.
Therefore, we define 
$\;{\cal R}\,:=\,\{\,\rho[R_z]\;|\;z\in{\Bbb R}\,\}\;$ 
and $\;{\cal S}\,:=\,\{\,\rho[S_z]\;|\;z\in{\Bbb R}\,\}\;$ 
where for every $\,z\in{\Bbb R}\,$ infinite sets $\,R_z,S_z\subset{\Bbb N}\,$
are defined so that 
if $\,x<y\,$ then on the one hand
$\;R_x\supset R_y\;$ and $\;R_x\setminus R_y\;$ is an infinite set,
and on the other hand $\;S_x\cap S_y\;$ is a finite set.
(For example, choose a bijection $\,\varphi\,$ 
from $\,{\Bbb N}\,$ onto $\,{\Bbb Q}\,$
and put $\;R_x\,:=\,\{\,n\in{\Bbb N}\;|\;x\leq \varphi(n)\,\}\;$
for every $\,x\in{\Bbb R}\,$.
Furthermore, for every $\,x\in{\Bbb R}\,$ choose 
a set $\;T_x\,\subset\,{\Bbb Q}\cap[x-1,x]\;$ with 
$\,T_x'=\{x\}\,$  and put $\,S_x:=\varphi^{-1}(T_x)\,$.)
Clearly, for every infinite set $\,S\subset{\Bbb N}\,$
the family $\;\{\,B\cup\{\,a_n\;|\;k\leq n\in S\,\}\;|\;
B\in{\cal B}_2\;\land\;k\in{\Bbb N}\,\}\;$ is 
a filter base which generates a filter $\,{\cal F}\,$
such that $\;{\cal N}_\eta(0)\cap{\cal F}\,=\,{\cal N}_{\rho[S]}(0)\,$.
Thus $\;{\cal R},{\cal S}\,\subset\,{\cal L}_0^*\;$ if $\,{\cal B}_2\,$ is 
countable. (So $\;{\cal R},{\cal S}\,\subset\,{\cal L}_0^*\;$ can be achieved 
if $\,\tau_2\in{\cal L}_0^*\,$.) If $\,\tau_2=\eta\,$ 
(and hence $\,{\cal B}_2=\{\emptyset\}\,$)
then the topologies in $\,{\cal R}\cup{\cal S}\,$ are homeomorphic.
Because if $\,S\subset{\Bbb N}\,$ is infinite 
then any 
increasing bijection from $\,{\Bbb R}\,$ onto $\,{\Bbb R}\,$
which maps $\,0\,$ to $\,0\,$ and 
$\;\{\,a_1,a_2,a_3,...\,\}\;$ onto $\;\{\,a_n\;|\;n\in S\,\}\;$
is clearly a homeomorphism from the space $\,({\Bbb R},\rho[{\Bbb N}])\,$
onto $\,({\Bbb R},\rho[S])\,$.
So in order to conclude the proof it 
remains to define a family $\,{\cal T}\,$ as desired.
\smallskip
For every free ultrafilter $\,{\cal F}\,$
on $\,{\Bbb N}\,$ put $\,\rho[{\cal F}]\,:=\,\bigcup_{S\in{\cal F}}\rho[S]\,$.
Clearly, $\;\tau_1\subset\rho[{\cal F}]\subset\tau_2\,$.
We claim that $\,\rho[{\cal F}]\,$ is a topology in the lattice $\,{\cal L}_0\,$.
Firstly, let  $\;U_1,U_2\in\rho[{\cal F}]\,$. Then
$\,U_i\in\rho[S_i]\,$ for $\,S_i\in{\cal F}\,$.
Since $\,S_1\cap S_2\,$ is an infinite set in the ultrafilter $\,{\cal F}\,$
and $\,\rho[S_1\cap S_2]\,$ is a topology containing
$\,\rho[S_1]\,$ and $\,\rho[S_2]\,$, the intersection $\,U_1\cap U_2\,$
lies in $\,\rho[S_1\cap S_2]\,$ and hence in 
$\,\rho[{\cal F}]\,$. Since $\;U\in\rho[S]\;$
whenever $\,0\not\in U\in\eta\;$ and $\,S\in{\cal F}\,$,
it is plain that the family $\,\rho[{\cal F}]\,$ is closed
under arbitrary unions and furthermore that $\,\rho[{\cal F}]\in{\cal L}_0\,$.
We also observe that for $\,U\in\eta\cap{\cal N}_\eta(0)\,$
we have $\,U\in\rho[{\cal F}]\,$ if and only if 
$\;U\supset B\;$ for some $\,B\in{\cal B}_2\,$
and $\;\{\,n\in{\Bbb N}\;|\;a_n\in U\,\}\in{\cal F}\,$.
Let $\,{\cal F}_1,{\cal F}_2\,$ be free ultrafilters 
on $\,{\Bbb N}\,$
and $\,S\in{\cal F}_1\,$
and assume $\;\rho[{\cal F}_1]\subset\rho[{\cal F}_2]\,$.
The set $\;V\,:=\,{\Bbb R}\setminus\{\,a_n\;|\;n\not\in S\,\}\;$
is $\tau_2$-open and $\;\{\,n\in{\Bbb N}\;|\;a_n\in V\,\}=S\,$.
Thus $\,V\,$ is $\rho[{\cal F}_1]$-open
and hence $\rho[{\cal F}_2]$-open and this implies 
$\,S\in{\cal F}_2\,$.
So we derive $\,{\cal F}_1\subset{\cal F}_2\,$ and hence
$\,{\cal F}_1={\cal F}_2\,$. 
Thus the topologies $\,\rho[{\cal F}]\,$ are mutually incomparable
and hence a family $\,{\cal T}\,$ as desired exists
provided that we always have $\,\rho[{\cal F}]\not\in{\cal L}_0^*\,$.
\smallskip
Assume for a contradiction that $\,\rho[{\cal F}]\in{\cal L}_0^*\,$
for a free ultrafilter $\,{\cal F}\,$ on $\,{\Bbb N}\,$.
Then we can choose a countable filter base 
$\;\{\,B_1,B_2,B_3,...\,\}\;$ of subsets of $\;{\Bbb R}\setminus[-1,1]\;$
such that $\;B_n\supset B_{n+1}\;$ for every $\,n\in{\Bbb N}\,$
and 
$\;\eta\cap{\cal N}_{\rho[{\cal F}]}(0)\,=\,
\{\,U\in\eta\;|\;\exists\,n\in{\Bbb N}:\;U\,\supset\,\{0\}\cup B_n\,\}\,$.
Put $\;S_m\,:=\,\{\,n\in{\Bbb N}\;|\;a_n\in B_m\,\}\;$ for every $\,m\in{\Bbb N}\,$.
Trivially, $\;S_{m}\supset S_{m+1}\;$ for every $\,m\in{\Bbb N}\,$.
Let $\,S\,$ be any set in the ultrafilter $\,{\cal F}\,$.
Then the set $\;{\Bbb R}\setminus\{\,a_n\;|\;n\not\in S\,\}\;$ is 
$\rho[{\cal F}]$-open and hence it contains $\,B_m\,$ for some $\,m\in{\Bbb N}\,$.
So for some $\,m\in{\Bbb N}\,$ we have 
$\;B_m\cap\{\,a_n\;|\;n\not\in S\,\}=\emptyset\;$ 
and hence $\;S_m\subset S\,$. Therefore, $\;\{\,S_m\;|\;m\in{\Bbb N}\,\}\;$
is a filter base for the filter $\,{\cal F}\,$. But this 
is impossible because a filter base for a free ultrafilter on $\,{\Bbb N}\,$
must be uncountable (cf.~[1] 7.8.a). 
This concludes the proof of Theorem 8.
\bigskip
{\it Remark.} For achieving $\,{\cal R},{\cal S}\subset{\cal L}_0^*\,$,
the additional 
assumption $\,\tau_2\in{\cal L}_0^*\,$ is essential in view of the following
counterexample $\,(\tau_1,\tau_2)\,$.
Consider the topologies $\,\tau_1:=\tau(\{0\})\,$ and 
$\,\tau_1':=\tau(]0,{1\over 2}[)\,$ in the Boolean algebra 
$\,{\cal D}\subset{\cal L}_0^*\,$. Let $\,\tau_2\,$ be the supremum
of $\,\{\tau_1,\tau_1'\}\,$ in the lattice $\,{\cal L}_0\,$. 
We observe that if $\,\tau_1\not=\tau\in{\cal L}_0\,$ and 
$\;\tau_1\subset\tau\subset\tau_2\;$ 
then $\,\tau\not\in{\cal L}_0^*\,$.
(Because for every $\,k\in{\Bbb N}\,$ and 
every sequence $\,(u_n)\,$
with $\,0<u_n\leq {1\over 2}\,$
the set $\;]{-1,1}[\,\cup\bigcup_{n=k}^\infty ]n,n+u_n[\;$
lies in $\,\tau\setminus\tau_1\,$.)
In particular, $\,\tau_2\not\in{\cal L}_0^*\,$. 
Furthermore, this counterexample demonstrates that
neither $\,{\cal D}\,$ nor 
$\,{\cal L}_0^*\,$ is a sublattice of $\,{\cal L}_0\,$.
\medskip
The minimum $\,\theta=\bigcap{\cal L}_0\,$
of the complete lattice $\,{\cal L}_0\,$ lies 
in $\,{\cal L}_0^*\,$. Thus by Theorem 8 and since 
it is clear that 
$\;{\cal L}_0\,=\,\{\,\tau\in{\cal L}\;|\;\tau\supset\theta\,\}\,$,
the topology $\,\theta\,$
has no immediate successor in the lattice $\,{\cal L}_0\,$ or in 
the partially ordered set $\,({\cal L},\subset)\,$.
On the other hand, the following proposition shows that 
the maximum $\,\eta=\bigcup{\cal L}_0\,$
of the lattice $\,{\cal L}_0\,$ has $\,2^c\,$
immediate predecessors in the lattice $\,{\cal L}_0\,$
which are also immediate predecessors of $\,\eta\,$
in the partially ordered family $\,({\cal L},\subset)\,$.
\medskip
{\bf Proposition 8.} {\it There exist $\,2^c\,$ (mutually non-homeomorphic)
topologies $\;\vartheta\in{\cal L}_0\;$ such that 
no topology from $\,{\cal L}\,$
lies strictly between $\,\vartheta\,$ and $\,\eta\,$.}
\smallskip
{\it Proof.} For a free ultrafilter $\,{\cal F}\,$
on $\,{\Bbb Z}\,$ let 
$\,\tau[{\cal F}]\,$ denote the topology as defined 
in the proof of Theorem~1. If
$\;{\cal K}\,\subset\,{\cal L}_0\setminus\{\eta\}\;$ 
is a chain with $\,\tau[{\cal F}]\in{\cal K}\,$
then $\,\eta\not=\bigcup{\cal K}\in{\cal L}_0\,$ 
by Proposition~6. Therefore, by applying Zorn's lemma,
for every free ultrafilter $\,{\cal F}\,$ on $\,{\Bbb Z}\,$ 
we can choose a maximal element $\,\vartheta[{\cal F}]\,$ in
the partially ordered set $\,({\cal L}_0\setminus\{\eta\},\subset)\,$
such that $\,\tau[{\cal F}]\subset \vartheta[{\cal F}]\,$.
For distinct free ultrafilters $\,{\cal F}_1,{\cal F}_2\,$
we have $\,\vartheta[{\cal F}_1]\not=\vartheta[{\cal F}_2]\,$
because $\,\tau[{\cal F}_1]\not=\tau[{\cal F}_2]\,$
and $\,\eta\,$ is the supremum 
of $\,\{\tau[{\cal F}_1],\tau[{\cal F}_2]\}\,$ in the lattice 
$\,{\cal L}_0\,$ in view of Proposition 1 since there are
sets $\,U_i\in\tau[{\cal F}_i]\,$
with $\;U_1\cap U_2\,=\;]{-1,1}[\,$.
(For example, choose $\;S_1\,\in\,{\cal F}_1\setminus{\cal F}_2\;$
and with $\;S_2\,:=\,{\Bbb Z}\setminus S_1\;$
put 
$\;U_i\,=\;]{-1,1}[\,\cup\,\bigcup_{n\in S_i}]n-{1\over 2},n+{1\over 2}[\;$
for $\,i\in\{1,2\}\,$.)
Finally, if $\,\eta\not=\tau\in{\cal L}\,$
and $\,\tau\supset\vartheta[{\cal F}]\,$ 
then $\,\tau=\vartheta[{\cal F}]\,$ 
since $\;{\cal L}_0\,=\,
\{\,\tau\in{\cal L}\;|\;\tau\supset\theta\,\}\;$
and $\,\vartheta[{\cal F}]\,$ is maximal in 
$\;{\cal L}_0\setminus\{\eta\}\,$, {\it q.e.d.}

\medskip\smallskip
{\it Remark.} By virtue of Theorem~8 every immediate predecessor
of $\,\eta\,$ in $\,{\cal L}_0\,$ must 
lie in $\,{\cal L}_0\setminus{\cal L}_0^*\,$.
This observation has two consequences in view of Proposition 8. Firstly we 
can be sure that 
$\;|{\cal L}_0\setminus{\cal L}_0^*|=|{\cal L}_0^*|=2^c\,$.
Secondly, the central assumption $\,\tau_1\in{\cal L}_0^*\,$ in 
Theorem~8 cannot
be replaced with the weaker assumption $\,\tau_1\in{\cal L}_0\,$.
\bigskip\smallskip
{\bff 13. Extremely long chains of topologies}
\medskip
Since both the 
existence of free ultrafilters
and the existence of the topologies $\,\vartheta[{\cal F}]\,$ 
in the proof of Proposition 8 are based 
on a maximality principle equivalent with the Axiom of Choice,
one might ask whether in the proof of Proposition 8
the topology $\,\tau[{\cal F}]\,$ is 
maximal in $\,{\cal L}_0\setminus\{\eta\}\,$ already, whence 
$\,\vartheta[{\cal F}]=\tau[{\cal F}]\,$.
This would be far from being true in view of the following theorem
which affirmatively answers the interesting question whether 
the lattice $\,{\cal L}_0\,$ contains chains of size 
greater than $\,c\,$. Define $\,\lambda:=\log(c^+)\,$,
i.e.~$\,\lambda\,$ is the smallest cardinal number $\,\kappa\,$
satisfying $\,2^\kappa>c\,$, whence $\,\aleph_1\leq\lambda\leq c\,$
and $\,2^\lambda>c\,$.
\medskip
{\bf Theorem 9.} %{\it Let $\,\lambda\,$ denote the smallest cardinal number 
%satisfying $\,2^\lambda>c\,$. (In standard notation, $\,\lambda=\log(c^+)\,$,
% whence $\,\aleph_1\leq\lambda\leq c\,$.)
{\it For every free ultrafilter $\,{\cal F}\,$ on 
$\,{\Bbb Z}\,$ there is a chain 
$\,{\cal K}\subset{\cal L}_0\,$ such that $\,|{\cal K}|=2^\lambda\,$
and $\;\tau\supset\tau[{\cal F}]\;$ for every $\,\tau\in{\cal K}\,$.}
\medskip
{\it Proof.}
For $\,n\in{\Bbb N}\,$ define a strictly increasing real function 
$\,\varphi_n\,$ by $\;\varphi_n(x)\,=\,3^{-n}(x+1)\,$, whence 
$\,\varphi_n\,$ maps $\,[0,1]\,$ onto $\;[3^{-n},2\cdot 3^{-n}]\,$.
For every set $\,A\subset[0,1]\,$ define 
\medskip
\centerline{$\Phi(A)\;\,:=\;\,\{0\}\;\cup\,\bigcup\limits_{k\in {\Bbb Z}}
\big(k+\!\bigcup\limits_{n=1}^\infty \varphi_n(A)\,\big)\,.$}
\smallskip
Let $\,{\cal F}\,$ be a free ultrafilter $\,{\cal F}\,$ on $\,{\Bbb Z}\,$.
For $\,A\subset[0,1]\,$
let $\,\tau[{\cal F},A]\,$ denote the coarsest topology in the lattice 
$\,{\cal L}_0\,$ which is finer than $\,\tau[{\cal F}]\,$ and contains
all Euclidean open sets 
$\;U\,\supset\,\Phi(A)\,$.
(In particular, $\,\tau[{\cal F},\emptyset]=\eta\,$.)
Since $\,\tau[{\cal F},A]\in{\cal L}_0\,$, it is plain 
that $\,W\in\eta\,$ is an open neighborhood of $\,0\,$ in the space
$\,({\Bbb R},\tau[{\cal F},A])\,$ if and only if 
$\;W=U\cap V\;$ for some $\,U,V\in\eta\,$ with $\,U\supset\Phi(A)\,$
and $\,0\in V\,$ and $\,V\cap{\Bbb Z}\in{\cal F}\,$.
\smallskip
Obviously, $\;\tau[{\cal F},B]\subset \tau[{\cal F},A]\;$
if $\;\emptyset\not=A\subset B\subset [0,1]\,$. 
Moreover, $\,\tau[{\cal F},B]\,$ is strictly coarser 
than $\,\tau[{\cal F},A]\,$ if 
$\;\emptyset\not=A\subset B\subset [0,1]\;$ and $\,A\not= B\,$. 
Because if $\;b\,\in\,B\setminus A\;$
then the Euclidean open set $\;Y\,:=\;]{-1,1}[\,\cup
({\Bbb R}\setminus({\Bbb Z}\cup\Phi(\{b\})))\;$
is $\tau[{\cal F},A]$-open since $\;Y\supset\Phi(A)\,$. But $\,Y\,$ is not 
$\tau[{\cal F},B]$-open because if $\,k\in{\Bbb Z}\,$ and $\,|k|\geq 2\,$
and $\,\varepsilon>0\,$ then $\,Y\,$ does not contain 
$\;]k,k+\varepsilon[\,\cap\,\Phi(B)\,$. 
\smallskip
Therefore, 
$\;{\cal K}\,=\,\{\,\tau[{\cal F},A]\;|\;A\in{\cal A}\,\}\;$
is a chain as desired if $\,{\cal A}\,$ 
is a chain of subsets of $\,[0,1]\,$
with $\,|{\cal A}|=2^\lambda\,$. 
\smallskip
Such a chain $\,{\cal A}\,$ can easily be defined as follows.
Choose a linearly ordered set $\,(L,\preceq)\,$
such that $\,|L|=2^\lambda\,$
and $\,L\,$ has a dense subset $\,D\,$ with $\,|D|=c\,$.
(This choice is possible in view of [1] Theorems 5.7.c and 5.8.b.)
Define a bijection $\,g\,$ from $\,D\,$ onto $\,[0,1]\,$
and put $\;A_x\,:=\,\{\,g(y)\;|\;x\prec y\in D\,\}\;$
for every $\,x\in L\,$. Finally define 
$\;{\cal A}\,:=\,\{\,A_x\;|\;x\in L\,\}\,$, {\it q.e.d.}
\bigskip
{\it Remark.} One does not need Theorem 9 to track down 
chains in $\,{\cal L}_0\,$ of size $\,2^\lambda\,$,
it is enough to define $\,{\cal A}\,$ as above
and to take into consideration 
 that our Boolean algebra $\,{\cal D}\subset{\cal L}_0^*\,$
is isomorphic with the power set of $\,[0,1]\,$.
The lattice $\,{\cal L}_0\,$ contains chains 
of the maximal possible size $\,2^c\,$ provided that $\,2^\lambda=2^c\,$.
Of course,  $\,2^\lambda=2^c\,$ trivially follows from 
the irrefutable hypothesis $\,\lambda=c\,$.
(Conversely, $\,2^\lambda=2^c\,$ does not imply $\,\lambda=c\,$.)
The hypothesis $\,\lambda=c\,$ is irrefutable
because $\,\lambda=c\,$ is obviously a 
consequence of the Continuum Hypothesis 
$\,\aleph_1=c\,$. However, the hypothesis $\,\lambda=c\,$ 
is much weaker than the very restrictive  
hypothesis $\,\aleph_1=c\,$
because it is consistent with ZFC set theory that $\,\lambda=c\,$ 
and $\,\aleph_1<\mu<c\,$ for infinitely many cardinal numbers $\,\mu\,$.
Even more, roughly speaking, $\,\lambda=c\,$ cannot prevent 
an arbitrarily large deviation of $\,c\,$ from $\,\aleph_1\,$.
(Precisely, in view of [3] 16.13 and 16.20, if 
$\,\kappa>\aleph_1\,$ is an arbitrary
regular cardinal in G\"odel's Universe
$\,{\rm L}\,$ then there is 
a generic extension $\,{\rm E}_\kappa\,$ of $\,{\rm L}\,$
preserving all cardinals such that $\,\lambda=c=\kappa\,$ holds
in the ZFC-model $\,{\rm E_\kappa}\,$.)
\bp
\hrule
\mp
Up to this point the present paper is essentially identical 
with the author's paper {\it Coarse topologies on the real line}
published in {\sl Matemati\v cki Vesnik {\bf 68}} (2016).
The next section is a supplement written in 2026. 
\mp
\hrule
\bigskip
{\bff 14. Counting with respect to the weight}
\mp
Referring to the remark in Section 4, 
all spaces constructed in the proof of Theorem 1 are of uncountable weight.
There arises the question whether they are of weight $\,c\,$.
Since for $\,\tau\in{\cal L}\,$ the weight of $\,(\R,\tau)\,$ 
cannot exceed $\,c\,$, the question can trivially be answered 
in the affirmative if the Continuum Hypothesis is assumed. 
In view of the following basic estimate (see [8] for a proof)
we can be sure that among the spaces depicted in Theorem~1 there 
actually are $\,2^c\,$ ones of weight 
$\,c\,$ provided that $\,2^\kappa<2^c\,$ for every cardinal number 
$\,\kappa<c\,$. As usual, $\,w(X)\,$ denotes the weight 
of the topological space $\,X\,$.
\mp
(14.1)$\;$ {\it If $\,\kappa\,$ is an infinite cardinal 
and $\,{\cal G}\,$ is a family of mutually non-homeomorphic
infinite Hausdorff spaces such that 
$\;\max\{|X|,w(X)\}\leq \kappa\;$
for every $\,X\in{\cal G}\,$ then $\;|{\cal G}|\leq 2^\kappa\,$.}
\mp
What can be accomplished if $\,2^\kappa=2^c\,$ for some $\,\kappa<c\,$?
In this case, can we rule out that all spaces 
depicted in Theorem 1 are of weight smaller than $\,c\,$
even when no space is second countable?
The most interesting question of course
is whether in (14.1) the maximal size $\,2^\kappa\,$ 
can be achieved for families $\,{\cal G}\,$  
of spaces as depicted in Theorem 1 and all possible 
weights $\,\kappa\,$ in view of the following fact.
\mp
(14.2)$\;$ {\it The existence of $\,c\,$ uncountable cardinals 
$\,\kappa<c\,$ is consistent with {\rm ZFC} set theory.} 
\mp
The consistency result (14.2) can easily be established
by routine forcing. (In G\"odel's universe L 
let $\,\theta\,$ be 
the smallest ordinal number 
of uncountable cofinality satisfying $\,\aleph_\theta=\theta\,$.
While $\,2^{\aleph_0}=\aleph_1\,$ holds in L,
there is a generic extension $\,V[{\rm L}]\,$ of L 
such that
$\,2^{\aleph_0}=\aleph_\theta\,$ holds in $\,V[{\rm L}]\,$, see [3] p.~226.) 
\mp
All the questions are answered by the following generalization 
of Theorem 1.
\mp
{\bf Theorem 10.}\quad {\it For every infinite cardinal $\,\kappa\leq c\,$
there exists a family $\;{\cal T}_\kappa\subset{\cal L}\;$ 
with $\,|{\cal T}_\kappa|=2^\kappa\,$ such that
$\,({\Bbb R},\tau)\,$ is a completely normal Baire 
space of weight $\,\kappa\,$ for each $\,\tau\in{\cal T}_\kappa\,$
and two spaces $\,({\Bbb R},\tau)\,$ and $\,({\Bbb R},\tau')\,$ are
never homeomorphic for distinct topologies 
$\,\tau,\tau'\in{\cal T}_\kappa\,$.}
\mp
In order to verify Theorem 10 we distinguish 
the cases $\,2^\kappa>c\,$ and $\,2^\kappa=c\,$.
(For no cardinal $\,\kappa\,$ with $\,\aleph_0<\kappa<c\,$ 
one can decide whether $\,2^\kappa>c\,$ or $\,2^\kappa=c\,$.)
The case $\,2^\kappa>c\,$ can be settled by a simple modification 
of the proof in Section 4. In the proof 
of Theorem~1 we worked with free ultrafilters   
for only two reasons, namely that a free ultrafilter on $\,\Z\,$
does not contain finite sets and that $\,\Z\,$ 
carries $\,2^c\,$ free ultrafilters, whence Lemma 3 can be applied.
The maximality of ultrafilters is also used to 
derive $\,\tau[{\cal F}_1]\not\subset\tau[{\cal F}_2]\,$
for distinct ultrafilters $\,{\cal F}_1,{\cal F}_2\,$.
Now, if we consider {\it free filters} on $\,\Z\,$
instead of free ultrafilters then we can still be sure that  
$\,\tau[{\cal F}_1]\not=\tau[{\cal F}_2]\,$
for distinct free filters $\,{\cal F}_1,{\cal F}_2\,$
and, naturally, that a free filter cannot contain finite sets.
Thus, by applying Lemma 3, the case $\,2^\kappa>c\,$ can be settled in 
view of the following proposition.
\mp
{\bf Proposition 9.} {\it For every infinite cardinal $\,\kappa\leq c\,$ 
there exist $\,2^\kappa\,$ free filters 
$\,{\cal F}\,$ on $\,\Z\,$ such that $\,\kappa\,$ is the 
least possible size of a filter base which generates $\,{\cal F}\,$.}
\mp\sp
For a filter $\,{\cal F}\,$ on an infinite 
set $\,S\,$ let $\,\chi({\cal F})\,$ denote the 
least possible size of a filter base which generates $\,{\cal F}\,$.
Trivially, $\;\chi({\cal F})\leq|{\cal F}|\leq 2^{|S|}\,$.
Now, Proposition 9 is identical with
the following statement for the special case $\,\kappa=\aleph_0\,$.
\medskip
(14.3) $\;$ {\it If $\,S\,$ is an infinite set and 
$\;|S|\leq \kappa\leq 2^{|S|}\;$ 
then there exist $\,2^{\kappa}\,$ free filters 
$\,{\cal F}\,$ on $\,S\,$ such that $\;\chi({\cal F})=\kappa\,$.}
\mp
For a proof of (14.3) see  [8] Prop.~3.
It is plain that the cardinal $\,2^{\kappa}\,$ in (14.3) 
and hence in Proposition 9 is best 
possible.
\mp
In order to finish the proof of Theorem 10 
we have to deal with the case $\,2^\kappa=c\,$.
In this case we cannot apply Lemma 3 and we actually have to 
track down $\,2^\kappa\,$ mutually non-homeomorphic spaces 
as desired. However, since $\,2^\kappa=c\,$, we have to track down 
{\it only} $\,c\,$ spaces as desired and hence, 
to accomplish this, it is natural 
to modify the proof of Theorem 4
by increasing the weight $\,\aleph_0\,$ of each 
real arc $\,A_H\,$ to $\,\kappa\,$ while keeping the arcs mutually 
non-homeomorphic. Well, the proof of Theorem 4 shows 
that for each $\,H\in {\cal H}\,$ we can fix a point 
$\,a_H\in A_H\,$ such that for distinct $\,H,H'\in{\cal H}\,$
not only $\,A_H\,$ and $\,A_{H'}\,$ but also 
the two spaces $\,A_H\setminus\{a_H\}\,$ and $\,A_{H'}\setminus\{a_{H'}\}\,$
are not homeomorphic. As a consequence, 
{\it there exists a family $\,{\cal L}'\subset{\cal L}\,$ 
such that $\,|{\cal L}'|=c\,$ and the space $\,(\R,\tau)\,$ 
is Polish for every $\,\tau\in{\cal L}'\,$ 
and for distinct $\,\tau_1,\tau_2\in{\cal L}'\,$ 
the subspace $\,\R\setminus\{0\}\,$ of $\,(\R,\tau_1)\,$ 
is not homeomorphic to the subspace $\,\R\setminus\{0\}\,$ of 
$\,(\R,\tau_2)\,$.} 
Now, let $\,{\cal F}\,$ be a free filter on $\,{\Bbb Z}\,$
with $\,\chi({\cal F})=\kappa\,$. 
For each $\,\tau\in{\cal L}'\,$ 
define a topology $\,\tau[{\cal F}]\,$ on $\,{\Bbb R}\,$
by declaring $\,U\subset{\Bbb R}\,$ open if and only if 
$\,U\,$ is open in $\,(\R,\tau)\,$ and satisfies $\;0\not\in U\;$
or $\;U\cap{\Bbb Z}\,\in\,{\cal F}\,$.
Then $\,\tau[{\cal F}]\,$ 
is a well-defined topology on $\,{\Bbb R}\,$ coarser than $\,\tau\,$
and hence coarser than $\,\eta\,$. 
By similar arguments as in Section 4, 
$\,({\Bbb R},\tau[{\cal F}])\,$ is a completely normal 
Baire space. Due to $\,\chi({\cal F})=\kappa\,$, 
the weight of the space $\,(\R,\tau[{\cal F}])\,$ is $\,\kappa\,$ 
or, equivalently, the least possible size of a neighborhood basis 
of the point $\,0\,$ in $\,(\R,\tau[{\cal F}])\,$ is $\,\kappa\,$.
On the other hand, the separable subspace $\,\R\setminus\{0\}\,$ 
of $\,(\R,\tau[{\cal F}])\,$ is metrizable and hence second countable.
Therefore, if $\,f\,$ is a homeomorphism from 
$\,(\R,\tau_1[{\cal F}])\,$ to 
$\,(\R,\tau_2[{\cal F}])\,$ for $\,\tau_1,\tau_2\in{\cal L}'\,$
then $\,f(0)=0\,$ and hence $\,f(\R\setminus\{0\})=\R\setminus\{0\}\,$
and hence $\,\tau_1=\tau_2\,$. This concludes the proof 
of Theorem 10.
Of course, the case $\,\kappa=\aleph_0\,$ (which is included in the 
case $\,2^\kappa=c\,$) is directly settled by Theorem 4. 
\mp\sp
{\it Remark.} Referring to the remark in Section 4, 
there are $\,2^c\,$ free ultrafilters $\,{\cal F}\,$ 
on $\,\Z\,$ such that the weight of $\,(\R,\tau[{\cal F}])\,$
is not only uncountable but actually equal to $\,c\,$. 
This is true because {\it on every infinite set $\,S\,$ 
precisely $\,2^{2^{|S|}}\,$ free ultrafilters 
$\,{\cal F}\,$ with $\,\chi({\cal F})=2^{|S|}\,$ exist.} 
(For an elementary proof see [8] Prop.~4.) 
\mp\sp
Of course, also Theorem 5 can be generalized as follows.
\mp
{\bf Theorem 11.}\quad {\it For every infinite cardinal number 
$\,\kappa\leq c\,$
there exist $\,2^\kappa\,$ 
mutually non-homeomorphic 
topologies $\,\tau\in{\cal L}\,$  
such that $\,({\Bbb R},\tau)\,$ is a completely normal
first category space of weight $\,\kappa\,$.}
\mp
The proof of Theorem 11 can be carried out similarly as the 
proof of Theorem 10. The case $\,2^\kappa>c\,$ in Theorem 11 is 
settled by modifying the proof of Theorem 5 in the same way as we 
modified the proof of Theorem 1 above. And to settle the case 
$\,2^\kappa=c\,$, in the same way as we 
modified the proof of Theorem 4 above, 
we implant one filter $\,{\cal F}\,$ with $\,\chi({\cal F})=\kappa\,$
in the constructions 
of the $\,c\,$ real arcs $\,h_H(\R)\,$ in the proof of Theorem 6 
in order to increase their weight 
to $\,\kappa\,$. 
\bigskip\bigskip
{\bff References}
\medskip
[1] Comfort, W.W., and Negrepontis, S.: {\it The Theory of Ultrafilters.}
Springer 1974.
\smallskip
[2] Engelking, R.: {\it General Topology}, revised and completed edition.
Heldermann 1989. 
\smallskip
[3] T.~Jech, {\it Set Theory}, 3rd ed. Springer 2002.
\smallskip
[4] Kuba, G.: {\it Counting metric spaces.}
Arch.~Math. {\bf 97},  569-578 (2011).
\smallskip
[5] Kuba, G.: {\it On certain 
separable and connected refinements of the Euclidean 
topology.} Matemati\v cki Vesnik {\bf 64} (2), 125-137 (2012).
\smallskip
[6] Kuba, G., {\it On the variety of Euclidean point sets.}
Internat.~Math.~News~{\bf 228}, 23-32 

\rightline{(2015) and {\tt arXiv:2004.11101v1}}
\smallskip
[7] Kuba, G., {\it Counting ultrametric spaces.} Coll.~Math.~{\bf 152},
217-234 (2018).
\sp
[8] Kuba, G.: {\it Counting overweight spaces.} {\tt arXiv:2006.02880v1}.
\smallskip
[9] Lelek, A., and McAuley, L.F.: {\it On hereditarily locally 
connected spaces and one-to-one continuous images of a line.}
Coll.~Math.~{\bf 17}, 319-324 (1967).
\bigskip\medskip
{\sl Author's address:} Institute of Mathematics, 
BOKU University, Vienna.
\smallskip
{\sl E-mail:} {\tt gerald.kuba(at)boku.ac.at}
\end